\documentclass{amsart}

\usepackage{amsmath, amsfonts, latexsym,amsthm, amsxtra, amssymb,graphicx,cite}

\def\NN{{\mathbb N}}
\def\ZZ{{\mathbb Z}}

\def\RR{{\mathbb R}}
\def\CC{{\mathbb C}}

\def\cH{{\mathcal H}} 
\def\cM{{\mathcal M}} 
\def\li{{\rm Li}_2} 
\def\cal{\mathcal}

\newtheorem{theorem}{Theorem}[section]

\newtheorem{lemma}[theorem]{Lemma}
\newtheorem{corollary}[theorem]{Corollary}
\newtheorem{proposition}[theorem]{Proposition}

\author{Changgui ZHANG}
\address{Laboratoire P. Painlev\'e  CNRS UMR 8524, UFR de Math\'ematiques,
Universit\'e des Sciences et Technologies de Lille, Cit\'e
scientifique, 59655 Villeneuve d'Ascq cedex, France.}
\email{zhang@math.univ-lille1.fr}
\title[A modular-type formula for $(x;q)_\infty$]{A modular-type formula for the infinite product $(1-x)(1-xq)(1-xq^2)(1-xq^3)\cdots$}
\date{October 3, 2011.
\textrm{{\bf MSC 2000 Subject Classifications:} 11F20, 11F27, 30D05, 33E05, 05A30.
{\bf Keywords:} Modular elliptic functions, Jacobi $\theta$-function, Dedekind $\eta$-function, Lambert series, $q$-series.
}}
\begin{document}
\maketitle

\begin{abstract}
Let $q=e^{2\pi i\tau}$, $\Im\tau>0$, $x=e^{2\pi i\xi}$, $\xi\in\CC$ and $(x;q)_\infty=\prod_{n\ge 0}(1-xq^n)$. Let $(q,x)\mapsto(q^*,\iota_q x)$ be the classical modular substitution given by the relations $q^*=e^{-2\pi i/\tau}$ and  $\iota_q x=e^{2\pi i\xi/{\tau}}$. The main goal of this paper is to give a modular-type representation for the infinite product $(x;q)_\infty$, this means, to compare the function defined by $(x;q)_\infty$ with that given by $(\iota_q x;q^*)_\infty$. Inspired by the work  \cite{St} of Stieltjes on semi-convergent series, we are led to a ``closed'' analytic formula for the ratio $(x;q)_\infty/(\iota_q x;q^*)_\infty$ by means of the dilogarithm combined with a Laplace type integral, which admits a divergent series  as Taylor expansion at $\log q=0$. Thus, the function $(x;q)_\infty$ is linked with its modular transform $(\iota_qx;q^*)_\infty$ in such an explicit manner that one can directly find the modular formulae known for Dedekind's eta function, Jacobi theta function, and also for certain Lambert series. Moreover, one can  remark  that our results allow Ramanujan's  formula \cite[Entry 6', p. 268]{Be1} (see also \cite[p. 284]{Ram}) to be completed as a convergent expression for the infinite product $(x;q)_\infty$.
\end{abstract}

\tableofcontents

\section{Introduction}

Let $q=e^{2\pi i\tau}$ and $\Im\tau>0$. For any  $x\in\CC$, let $(x;q)_\infty=\prod_{n\ge 0}(1-xq^n)$. This series may be found in Euler \cite[Chap. XVI] {Eu}, from that one knows the following identity:
\begin{equation}
(x;q)_\infty=\sum_{n\ge 0}\frac{q^{n(n-1)/2}}{(q;q)_n}\,(-x)^n;
\end{equation}
here and in the following, one writes
$(q;q)_0=1$ and $(q;q)_n=(1-q)...(1-q^n)$ for $n\ge 1$. It is in \cite{Ja} that this $q$-series is used by Jackson for the definition of the basic Gamma function $\Gamma_q(z)$, usually called the Jackson's Gamma function: for all $z\in\CC$ verifying $q^z\notin\{1,q^{-1},q^{-2},...\}$, one defines
\begin{equation}
 \Gamma_q(z)=\frac{(q;q)_\infty}{(q^z;q)_\infty}\,(1-q)^{1-z}.
 \end{equation}

As mentioned in \cite{Ja}, the infinite product $(x;q)_\infty$  already appeared, often in connection with elliptic function theory, in works of Weierstrass, Halphen, Tannery, {Heine}, Rogers, Barnes, {\it etc}. In almost the same manner as what done by the usual Gamma function in the classical theory of special functions, the Jackson's Gamma function plays an important role for basic hypergeometric function theory; see \cite{GR}. We content also to emphasize the book \cite[Chapiter 10]{AAR} in which one can find  references about the Gauss' $q$-binomial theorem, the convergence of $q$-Gamma function towards Euler Gamma function, the triple-product formula of Jacobi, among many other matters in relation with $(x;q)_\infty$.

In his {\it Notebooks} \cite{Ram} at the page 284, Ramanujan gave a formula that is  recoded as Entry 6' in \cite[p. 268]{Be1} in the following manner:
\begin{equation}\label{equation:Ram}
(x;q)_\infty=\frac{\sqrt{2\pi s(1-x)}}{\Gamma(s+1)}\,q^{-1/24}\,e^{s(\log s-1)}\,e^{\frac{\li(x)}{\log q}-\theta}\,,
\end{equation}
where  $s=\log x/\log q$, $\li$ denotes the dilogarithm and where
\begin{equation}\label{equation:Ramtheta}
\theta=\sum_{n=1}^\infty\frac{B_{2n}}{(2n)!}\,(\log q)^{2n-1}\,\bigl\{\frac{B_{2n}}{2n}\,\frac{\log x}{1!}+\frac{B_{2n+2}}{2n+2}\,\frac{(\log x)^3}{3!}+...\bigr\}\,.
\end{equation}
In the above, to make the notations consistent with the rest of the present paper, we replaced in \cite[(6.8)]{Be1} the expressions $e^{-a}$ and $e^{-x}$ with $x$ and $q$, respectively; at the same time, according to \cite[p. 284]{Ram}, we restored the exponent $(2n-1)$, instead of $(2n+1)$ as appeared in \cite[(6.8)]{Be1}, for $\log q$ in the expression of $\theta$. It is worth noting that, in certain strict sense,  the equality \eqref{equation:Ram} remains {\it incomprehensible}, since the left hand side represents clearly an analytic function for $\vert q\vert<1$ and $x\in\CC$ while the convergence of the power series expansion $\theta$ in the right hand side can not be assumed.

In our opinion, the remainder term $\theta$ of Ramnujan, if it would be analytic, contains not only the semi-convergent power series in $\log q=0$ as given in his text, but also a complementary term that is exponentially small with respect to $\log q$ for $q\to 1^-$. In other words, in order to make the formula \eqref{equation:Ram} completed analytically, one needs to replace $\theta$ with the Borel-sum of a semi-convergent series plus a flat function in the analytical scale for $\log q$. To do this, we shall consider both the function $(x;q)_\infty$ and its modular counterpart $(\iota_q x;q^*)_\infty$ with  $x=e^{2\pi i\xi}$, $q^*=e^{-2\pi i/\tau}$ and  $\iota_q x=e^{2\pi i\xi/{\tau}}$. With these notations, if $\alpha\in]0,1[$ or $\alpha=0$, the limit relation $q\to e^{2\alpha \pi i}$ is equivalent to say that $q^*\to e^{-2\pi i/\alpha}$ or $q^*\to 0$, respectively; in particular, as $q$ tends towards $1$ inside the unity circle, the modular variable $q^*$ is exponentially small.

In the following, we will begin with some classic modular relations known on $\eta$ and $\theta$-functions with together a work of Stieltjes on semi-convergent series; these old-fashion works constitute our first motivation for the present work. The plan of the rest of the paper will be outlined in \S\ref{subsection:plan}. To conclude this introduction, we will give in \S\ref{subsection:qdiff} some general commentaries for the analytic theory of $q$-difference equations with connection to the study of $q$-series and, in particular, to  some Ramanujan's dream.

\subsection{Revisit on  $\eta$, $\theta$-modular relations and semi-convergent series}The strategy that we shall put in place in the present work is initially inspired by the following observations.

Firstly, the infinite product $(q;q)_\infty$, intimately associated to the Dedekind's $\eta$ function, is classically known to satisfy a modular relation \cite[(44), p. 154]{Se} and this is just the value of the function $(x;q)_\infty$ taken at $x=q$. Secondly, the Jacobi $\theta$ function $\sum_{n\in\ZZ}q^{n^2/2}(-x)^n$ satisfies also a modular relation and this can be written as the product of  $(q;q)_\infty$ by the factor $(\sqrt qx;q)_\infty(\sqrt q/x;q)_\infty$, which is left invariant by the substitution $x\mapsto1/x$.

Furthermore, in the last paragraph of his Ph.D Thesis \cite[p. 252-258]{St} on semi-convergent series, in connection with the later theory of Borel-summable series \cite{Ra}, Stieltjes found an singular integral representation for a Lambert series and this expression could be seen as a modular formula for the studied series.

So, we are led to consider the following question. Could one make use of Stieltjes approach to obtain an explicit formula that relates the function $(x;q)_\infty$ and its {\it modular} counterpart $(\iota_q x;q^*)_\infty$, where  $x=e^{2\pi i\xi}$, $q^*=e^{-2\pi i/\tau}$, and  $\iota_q x=e^{2\pi i\xi/{\tau}}$~? More precisely, the expected formula would be expressed in such an explicit manner that one might immediately deduce from that the known modular relations for $\eta$ and $\theta$ as recalled in the above.

In this paper, we shall show that a such formula exists and that, up to an explicit factor, the function defined by the product $(x;q)_\infty$ can be seen somewhat modular. This {\it non-modular part} will be  represented by the above-mentioned expansion $\theta$ of Ramanujan in \eqref{equation:Ram}, which is a divergent but Borel-summable or called semi-convergent power series on variable $\log q$ near zero. These results, subjects of Theorems \ref{theorem:main} and \ref{theorem:xqmain}, give rise to one new and unified approach to treat Jacobi theta function, Lambert series and other $q$-series or $q$-functions such as the Jackson's Gamma function.

Remember finally the first non-trivial example of $q$-series may certainly be the infinite product $(q;q)_\infty$, that is considered in Euler \cite[Chap. XVI] {Eu} and then is revisited by many of his successors, particularly intensively by Hardy and Ramanujan \cite[p. 238-241; p. 276-309; p. 310-321]{Ha} for the theory of partition. Concerning  the Hardy-Ramanujan's formula on $p(n)$, which finally becomes completed by Rademacher (and Selberg), on can see the beautiful paper \cite{Sel}. In a similar manner as what happened for this formula $p(n)$, our modular-type formula \eqref{equation:xqmain} permits to give an analytic sense to the above-recalled relation \eqref{equation:Ram} of Ramanujan; see \S\ref{subsection:Ramanujan} and \S\ref{subsec:limit} in the below for more details.

\subsection{Plan of the paper}\label{subsection:plan}
The paper is organized as follows.  Section \ref{section:main} is devoted to \emph{sight-read} certain terms contained in Theorem \ref{theorem:main} of \S \ref{subsec:main}. Firstly, in \S \ref{subsec:mainbister}, we will give two equivalent formulations of Theorem \ref{theorem:main}, one of which will be used in complex plane in \S \ref{subsec:complex}. In \S \ref{subsec:M} and \S \ref{subsec:P}, we will observe that the modular relation remains almost valid but a perturbation term exists. In \S \ref{subsec:Stirling}, we deal with the remainder term of the Stirling asymptotic formula for $\Gamma$-function.

Theorem \ref{theorem:xqmain}, given in \S \ref{subsec:complex}, is another equivalent version of Theorem \ref{theorem:main} and will be used in Section \ref{section:modular}, as it is formulated in terms of complex variables. Relation \eqref{equation:doubleGamma} shows that the above-mentioned non-modular part can be expressed in terms of the quotient of two Barnes' double Gamma functions. Finally, we will give, in Theorem \ref{theorem:Ram}, a complete analytic version for Ramanujan's formula \eqref{equation:Ram}, in which the formal power series $\theta$ will be represented by a function analytic on some open sector. In \S \ref{subsec:limit}, we will give some remarks about the limit behavior when $q$ goes to one by real values.

In Section \ref{section:modular}, we will explain how to utilize Theorem \ref{theorem:main} to get the classical modular formula for eta or theta function. In \S \ref{subsec:eta}, it will be merely observed that the so-called non-modular part identically vanishes; in the theta function case (\S \ref{subsec:theta}), two non-modular parts are of opposite sign  and then cancel each other out. In \S \ref{subsec:anotherproof}, a second proof will be delivered to $\theta$-modular equation from the point view of $q$-difference equations. In \S \ref{subsec:L12}, we consider the first order derivatives of $(1-x)(1-xq)(1-xq^2)(1-xq^3)...$ and then get some results for two families of $q$-series, including Lambert series as special cases that will be treated in \S \ref{subsec:Lambert}.

In Section \ref{section:proof}, we give a complete proof of our main Theorem. To do this, we need several elementary but somewhat technical calculations, that will be formulated in terms of various lemmas.

\subsection{Analytic theory of $q$-difference equations and Ramanujan's dream}\label{subsection:qdiff}
We are interested in studying the analytical theory of differential, difference and $q$-difference equations, \emph{\`a la Birkhoff} \cite{Bi}; see \cite{DZ}, \cite{RSZ}, \cite{zh}. The elliptic functions and one variable modular functions can be seen as specific solutions of certain particular difference  or $q$-difference equations; in this line, we shall give a proof on Theta function modular equation in \S \ref{subsec:anotherproof}. We believe that a good understanding of singularities structure, that is, Stokes analysis  \cite{Ra} as well as other geometric tools, often permits a lot more of comprehension about certain magical formulas or, say, some Ramanujan's dream.

The main results of this paper are announced in \cite{zh4}.

\section{Modular-type expansion of $(x;q)_\infty$}\label{section:main}

Let $q$ and $x$ be complex numbers; if $|q|<1$, we recall that
$$
(x;q)_\infty=\prod_{n=0}^\infty(1-xq^n)\,.
$$

In \S \ref{subsec:main} and \S \ref{subsec:mainbister}, we will suppose that $q\in(0,1)$ and $x\in(0,1)$, so that the infinite product $(x;q)_\infty$ converges in $(0,1)$; therefore, one can take the logarithm of this function. From \S \ref{subsec:M}, we will work with complex variables.
In \S \ref{subsec:complex},  a modular-type expansion for $(x;q)_\infty$ will be given in complex plane.
As usual, $\log$ will stand for the principal branch of the logarithmic function over its Riemann surface denoted by $\tilde\CC^*$ and, in the meantime, the broken plane $\CC\setminus(-\infty,0]$ will be identified to a part of $\tilde\CC^*$.

\subsection{Statement of one main result in $\RR$}\label{subsec:main}

The main result of our paper is the following statement.

\begin{theorem}\label{theorem:main}
Let $q=e^{-2\pi\alpha}$, $x=e^{-2\pi(1+\xi)\alpha}$ and suppose $\alpha>0$ and $\xi>-1$. The following relation holds:
\begin{align}\label{equation:main}
\log(x;q)_\infty=&-\frac{\pi}{12\alpha}+\log\frac{\sqrt{2\pi}}{\Gamma(\xi+1)}+\frac \pi{12}\,\alpha-\bigl(\xi+\frac12\bigr)\log\frac{1-e^{-2\pi\xi \alpha}}\xi\cr
&
+\int_0^{\xi }\bigl(\frac {2\pi\alpha t}{e^{2\pi\alpha t}-1}-1\bigr)\,dt+M(\alpha,\xi),
\end{align}
where
\begin{equation}\label{equation:M}
M(\alpha,\xi)=-\sum_{n=1}^\infty\frac{\cos2n\pi\xi}{n(e^{2n\pi/\alpha}-1)}-\frac2\pi\,{\cal PV}\int_0^\infty\sum_{n=1}^\infty\frac{\sin2n\xi \pi t}{n(e^{2n\pi t/\alpha}-1)}\,\frac{dt}{1-t^2} \,.
\end{equation}
\end{theorem}

In the above, $\displaystyle {\cal PV}\int$ stands for the principal value of a singular integral in the Cauchy's sense; see \cite[\S 6.23, p. 117]{WW} or the corresponding definition recalled later in \S \ref{subsec:RN23}. We will leave the proof  until Section \ref{section:proof}.

Before extending the main theorem to the complex plane (\S \ref{subsec:complex}), we first give some equivalent statements.

\subsection{Variants of Theorem \ref{theorem:main}}\label{subsec:mainbister}
Throughout all the paper,  we let
\begin{equation}\label{equation:B}
B(t)=\frac1{e^{2\pi t}-1}-\frac1{2\pi t}+\frac12\,;
\end{equation}
therefore, Theorem \ref{theorem:main} can be stated as follows.

\begin{theorem}\label{theorem:mainbis}
Let $q$, $x$, $\alpha$, $\xi$ and $M(\alpha,\xi)$ be as given in Theorem \ref{theorem:main}. Then the following relation holds:
\begin{align}\label{equation:mainbis}
\log(x;q)_\infty=&-\frac{\pi}{12\alpha}-\bigl(\xi+\frac12\bigr)\,\log2\pi\alpha+\log\frac{\sqrt{2\pi}}{\Gamma(\xi+1)}+\frac\pi2\,(\xi+1)\,\xi\alpha\cr
&+\frac \pi{12}\,\alpha+2\pi\alpha\int_0^\xi\bigl(t-\xi-\frac12\bigr)\,B(\alpha t)\,dt +M(\alpha,\xi),
\end{align}
where $B$ is defined in \eqref{equation:B}.
\end{theorem}

\begin{proof}
It suffices to notice the following elementary integral: for any real numbers $\lambda$ and $\mu$,
\begin{equation}\label{equation:integralB}
\int_0^\lambda B(\mu t)dt=\frac1{2\pi\mu}\,\log\frac{1-e^{-2\pi \lambda\mu}}{2\pi \lambda\mu }+\frac \lambda2\,.
\end{equation}

\end{proof}

As usual, let $\li$ denote the dilogarithm function; recall $\li$ can be defined as follows \cite[(2.6.1-2), p. 102]{AAR}:
\begin{equation}\label{equation:dilog}
\li (x)=-\int_0^x\log(1- t)\,\frac{dt}{t}=\sum_{n= 0}^{\infty}\frac{x^{n+1}}{(n+1)^2}\,.
\end{equation}

\begin{theorem}\label{theorem:mainter}
The following relation holds for any $q\in(0,1)$ and $x\in(0,1)$:
\begin{align}\label{equation:mainter}
\log(x;q)_\infty=&\frac1{\log q}\,{\li(x)}+\log\sqrt{1-x}-\frac{\log q}{24}\cr
&-\int_0^\infty B(-\frac{\log q}{2\pi}\,t)\,x^{t}\,\frac{dt}t+M(-\frac{\log q}{2\pi}\,,\log_qx)\,,
\end{align}
where $B$ denotes the function given by \eqref{equation:B}.
\end{theorem}

\begin{proof}
By the first Binet integral representation stated in \cite[Theorem 1.6.3 (i), p. 28]{AAR} for $\log\Gamma$, Theorem \ref{theorem:main} can be formulated as follows:
\begin{align}\label{equation:mainterI}
\log(x;q)_\infty=&-\frac{\pi}{12\,\alpha}\,+\frac\pi{12}\,\alpha-I_\Gamma(\xi)-\frac12\,\log(1-e^{-2\pi\xi\alpha}\bigr)\cr
&-\int_0^\xi\log(1-e^{-2\pi\alpha t})\,dt+M(\alpha,\xi),
\end{align}
where $\xi=\log_q(x/q)$ and $I_\Gamma$ denotes the corresponding Binet integral in term of the function $B$ defined by \eqref{equation:B}:
\begin{equation}\label{equation:IGamma}
I_\Gamma(\xi)=\log\Gamma(\xi+1)-\bigl(\xi+\frac12\bigr)\,\log\xi+\xi-\log\sqrt{2\pi}=\int_0^\infty B(t)\,e^{-2\pi\,\xi t}\, \frac{dt}{t}\,.
\end{equation}
Write $\log(x/q;q)_\infty=\log(1-x/q)+\log(x;q)_\infty$, and substitute $q$ by $e^{-2\pi \alpha}$ and $x/q=e^{-2\pi\xi\alpha}=q^\xi$ by $x$ in \eqref{equation:mainterI}, respectively;  we arrive at once at the following expression:
\begin{align*}
\log(x;q)_\infty=&\frac{\pi^2}{6\,\log q}-\frac{\log q}{24}+\log\sqrt{1-x}-\int_0^{\log_qx}\log\bigl(1-q^t\bigr)\,dt\cr
&-\int_0^\infty B(-\frac{\log q}{2\pi}\,t)\,x^{t}\,\frac{dt}t+M(-\frac{\log q}{2\pi}\,,\log_qx)\,,
\end{align*}
from which, using \eqref{equation:dilog}, we easily deduce the expected formula \eqref{equation:mainter}, for $\li(1)=\displaystyle\frac{\pi^2}6$.
\end{proof}

\subsection{Almost modular term $M$}\label{subsec:M}

We shall write the singular integral part in \eqref{equation:M} by means of contour integration in the complex plane, as explained in \cite[\S 6.23, p. 117]{WW}.
Fix a real $r\in(0,1)$ and let $\ell_r^-$ (resp. $\ell_r^+$) denote the path that goes along the positive axis from the origin $t=0$ to infinity via the half circle starting from $t=1-r$ to $1+r$ {\it below} (resp. {\it over}) its center point $t=1$. Define $P^\mp(\alpha, \xi)$ as follows:
\begin{equation}\label{equation:Pmp}
P^\mp(\alpha,\xi):=-\frac2\pi\,\int_{\ell_r^\mp}\sum_{n=1}^\infty\frac{\sin2n\xi \pi t}{n(e^{2n\pi t/\alpha}-1)}\,\frac{dt}{1-t^2}\,,
\end{equation}
where $\alpha>0$ and where $\xi$ may be an arbitrary  real number.

Observe that the integral on the right hand side of \eqref{equation:Pmp} is independent of the choice of $r\in(0,1)$, so that we leave out the parameter $r$ from  $P^\mp$.  Moreover, the principal value of the singular integral considered in \eqref{equation:M} is merely the average of $P^+$ and $P^-$, that is to say:
\begin{equation}\label{equation:sumP}
P^-(\alpha,\xi)+P^+(\alpha,\xi)=-\frac4\pi\,{\cal PV}\int_0^\infty\sum_{n=1}^\infty\frac{\sin2n\xi \pi t}{n(e^{2n\pi t/\alpha}-1)}\,\frac{dt}{1-t^2}\,.
\end{equation}
By the residues Theorem, we find:
\begin{equation}\label{equation:StokesP}
P^-(\alpha,\xi)-P^+(\alpha,\xi)={2i}\,\sum_{n=1}^\infty\frac{\sin2n\xi \pi}{n(e^{2n\pi /\alpha}-1)}\,,
\end{equation}
from which we arrive at the following expression:
\begin{equation}\label{equation:MP0}
M(\alpha,\xi)=P^-(\alpha,\xi)-\sum_{n=1}^\infty\frac{e^{2n \pi\xi i}}{n(e^{2n\pi /\alpha}-1)}\,.
\end{equation}

\begin{theorem}\label{theorem:MP}
Let $M$ be as in Theorem \ref{theorem:main} and let $P^-$ be as in \eqref{equation:Pmp}. For any $\xi\in\RR$ and $\alpha>0$, the following relation holds:
\begin{equation}\label{equation:MP}
M(\alpha,\xi)=\log(e^{2 \pi\xi i-2\pi /\alpha};e^{-2\pi /\alpha})_\infty+P^-(\alpha,\xi)\,,
\end{equation}
where $\log$ denotes the principal branch of the logarithmic function over its Riemann surface.
\end{theorem}

\begin{proof}
Relation \eqref{equation:MP} follows directly from \eqref{equation:MP0}. Indeed, for the last series of \eqref{equation:MP0}, one can expand each fraction $(e^{2n\pi /\alpha}-1)^{-1}$ as power series in $e^{-2n\pi /\alpha}$ and then permute the summation order inside the obtained double series, due to absolute convergence.
\end{proof}

Consequently, the term $M$ appearing in Theorem \ref{theorem:main} can be considered as being an \emph{almost modular term} of $\log(x;q)_\infty$;  the correction term $P^-$ given by \eqref{equation:MP} will be called {\it disruptive factor or perturbation term}.

\subsection{Perturbation term $P$}\label{subsec:P}

In view of the classical relation
\begin{equation}\label{equation:cot}
\cot\frac{t}2=\frac2{t}-\sum_{n=1}^\infty\frac{4t}{4\pi^2n^2-t^2}\,,
\end{equation}
from \eqref{equation:Pmp} one can obtain the following expression:
\begin{equation}\label{equation:P-}
P^-(\alpha,\xi)=\int_{\ell_r^{-*}}\frac{\sin\xi t}{e^{ t/\alpha}-1}\,\big(\cot\frac{t}2-\frac2{t}\bigr)\,\frac{dt}{t}\,.
\end{equation}
In the last integral \eqref{equation:P-},  $r\in(0,1)$ and
$$\ell^{-*}_r=(0,1-r)\cup\Bigl(\cup_{n\ge 1}\bigl(C_{n,r}^-\cup(n+r,n+1-r)\bigr)\Bigr)\,,
$$
 where for any positive integer $n$, $C_{n,r}^-$ denotes the half circle passing from $n-r$ to $n+r$ by the right hand side.

One may replace the integration path $\ell^{-*}_r$ by any half line from origin to infinity which does not meet the real axis. In view of what follows in matter of complex extension considered in \S \ref{subsec:complex}, let us first introduce the following modified complex version of $P^-$: for any $d\in(-\pi,0)$, let
\begin{equation}\label{equation:Pd}
P^d(\tau,\nu)=\int_{0}^{\infty e^{id}}\frac{\sin\frac\nu\tau\, t}{e^{i t/\tau}-1}\,\big(\cot\frac{t}2-\frac2{t}\bigr)\,\frac{dt}{t}\,,
\end{equation}
the path of integration being the half line starting from origin to infinity with argument $d$.

From then on, if we let $\tilde\CC^*$ to denote the Riemann surface of the logarithm, we will define
\begin{equation}\label{equation:Sab}
S(a,b):=\{z\in\tilde\CC^*:\arg z\in(a,b)\}
\end{equation}
 for any pair of real numbers $a<b$; notice that the Poincar\'e's half-plane $\cH$ will be identified to $S(0,\pi)$ while the  broken plane $\CC\setminus(-\infty,0]$ will be seen as the subset $S(-\pi,\pi)\subset\tilde\CC^*$.

\begin{lemma}\label{lemma:Pd}
The family of functions $\{P^d\}_{d\in(-\pi,0)}$ given by \eqref{equation:Pd} gives rise to an analytical function over the domain \begin{equation}\label{equation:Omega-}
\Omega_-:=S(-\pi\,,\pi)\times\Bigl( \CC\setminus\bigl((-\infty,-1]\cup[1,\infty)\bigr)\Bigr)\subset\CC^2\,.
\end{equation}
Moreover, if we denote this function by $P_-(\tau,\nu)$, then the following relation holds for all $\alpha>0$ and $\xi\in\RR$:
\begin{equation}\label{equation:PP-}
 P_-(\alpha i,\xi\alpha i)=P^{-}(\alpha,\xi)\,.
\end{equation}

\end{lemma}

\begin{proof}
Let $B$ be as in \eqref{equation:B}; from the relation
\begin{equation}\label{equation:cotB}
\cot\frac t2-\frac2t={2i}\,B(\frac{it}{2\pi})\,,
\end{equation}
it follows that the function $P^d$ given by \eqref{equation:Pd} is well defined and analytic at $(\tau,\nu)=(\tau_0,\nu_0)$ whenever the corresponding integral converges absolutely, that is, when the following condition is satisfied:
$$
\bigl\vert\Re(\frac{e^{id}}{\tau_0}\,\nu_0 i)\bigr\vert<\Re(\frac{e^{id}}{\tau_0}\,i)\,.
$$
Therefore, $P^d$ is analytic over the domain $\Omega^d$ if we set
\begin{align}\label{equation:Omegad}
\Omega^d=\cup_{\sigma\in(0,\pi)}\bigl(0,\infty e^{i(d+\sigma)}\bigr)\times\{\nu\in\CC:\bigl\vert\Im(\nu\,e^{-i\sigma})\bigr\vert<\sin\sigma\}\,.\end{align}

Thus we get the analyticity domain $\Omega_-$ and also relation \eqref{equation:PP-} by using the standard argument of analytic continuation.
\end{proof}

Let us give some precision about the above-employed continuation procedure, which is really a radial continuation. In fact, for any pair of directions of arguments $d_1$, $d_2\in(-\pi,0)$, say $d_1<d_2$, the common domain $\Omega^{d_1}\cap \Omega^{d_2}$ contains a (product) disk $D(\tau_0;r)\times D(0;r)$ for certain $\tau_0\in S(d_2,d_1+\pi)$ and some radius $r>0$, and all points in both $\Omega^{d_1}$ and $\Omega^{d_2}$ can be almost radially joined to this disk.

On the other hand, if we take the arguments $d\in(0,\pi)$ instead of $d\in(-\pi,0)$ in \eqref{equation:Pd}, we can get an analytical function, say $P_+$, defined over
$$
\Omega_+:=S(0\,,2\pi)\times\Bigl( \CC\setminus\bigl((-\infty,-1]\cup[1,\infty)\bigr)\Bigr)
$$ and such that, for all $\alpha>0$ and $\xi\in\RR$:
\begin{equation}\label{equation:PP+}
 P_+(\alpha i,\xi\alpha i)=P^{+}(\alpha,\xi)\,.
\end{equation}
Therefore, the Stokes relation \eqref{equation:StokesP} can be extended in the following manner.

\begin{theorem}\label{theorem:StokesP}
For any $\tau\in\cH$, the relation
\begin{equation}\label{equation:StokesP+-}
P_-(\tau ,\nu)-P_+(\tau ,\nu)={2i}\,\sum_{n=1}^\infty\frac{\sin\frac{2n\nu\pi}\tau}{n(e^{2n\pi i/\tau}-1)}
\end{equation}
holds provided that $|\Im(\nu/\tau)|<-\Im(1/\tau)$.
\end{theorem}

\begin{proof}
In view of \eqref{equation:PP-} and \eqref{equation:PP+}, one may observe that the expected relation \eqref{equation:StokesP+-} really reduces to \eqref{equation:StokesP} when $\tau=\alpha i$, $\nu=\xi i$, $\alpha>0$ and $\xi\in\RR$. Thus one can get \eqref{equation:StokesP+-} by an analytical continuation argument. Another way to arrive at the result is to directly use the residues theorem.
\end{proof}

Using \eqref{equation:cotB}, one can write \eqref{equation:Pd} as follows:
$$
P^d(\tau,\nu)=2i\,\int_{0}^{\infty e^{id}}\bigl(B(\frac{t}\tau\,i)-\frac {\tau i}{2\pi t}-\frac12\bigr)\,B(it)\,\sin\frac{2\pi\nu t}\tau\,\frac{dt}{t}\,,
$$
where $B$ denotes the odd function given by \eqref{equation:B}. We guess that this expression contains some \emph{modular} information about the perturbation term !

\subsection{Remainder term relating Stirling's formula}\label{subsec:Stirling}

Let us consider the integral term involving the function $B$ in formula \eqref{equation:mainter} of Theorem \ref{theorem:mainter}, which is, up to the sign, the remainder term $I_\Gamma$ appearing in the Stirling's formula; see \eqref{equation:IGamma}. So, we introduce the following family of associated functions: for any $d\not=\frac\pi 2 \bmod\pi$, define
\begin{equation}\label{equation:g}
 g^d(z)=-\int_0^{\infty e^{id}} B(t)\,e^{-2\pi z t}\,\frac{dt}{t}\,.
\end{equation}
It is obvious that $g^d$ is analytic over the half plane $S(-\frac\pi2-d,\frac\pi2-d)$, where $S(a,b)$ is in the sense of \eqref{equation:Sab}.  By usual analytic continuation, each of the families of functions $\{g^d\}_{d\in(-\frac\pi2,\frac\pi2)}$ and $\{g^d\}_{d\in(\frac\pi2,\frac{3\pi}2)}$ will give rise to a function that we denote by $g^+$ and $g^-$ respectively; that is, $g^+$ is defined and analytical over the domain $S(-\pi,\pi)$  while $g^-$,  over $S(-2\pi,0)$. Since $B(-t)=-B(t)$, it follows that
\begin{equation}\label{equation:g+-}
	g^+(z)=-g^-(e^{-\pi i}\,z)
	\end{equation}
for any $z\in S(-\pi,\pi)$. Moreover, if $z\in S(-\pi,0)$, one can choose a small $\epsilon>0$ such that $g^\pm(z)=g^d(z)$, $d={\pi/2\mp\epsilon}$; by applying the residues theorem to the following contour integral
$$
\bigl(\int_0^{\infty e^{i({\pi/2+\epsilon})}} -\int_0^{\infty e^{i({\pi/2-\epsilon})}}\bigr)\,B(t)\,e^{-2\pi z t}\,\frac{dt}{t}\,,
$$
we find:
\begin{align}\label{equation:g-+}
g^+(z)-g^-(z)=-2\pi i\,\sum_{n\ge 1}\frac{e^{-2\pi z(ni)}}{2\pi\,ni}=\log(1-e^{-2\pi iz})\,.
\end{align}

\begin{lemma}\label{lemma:g}
The following relations hold: for any $z\in S(-\pi,0)$,
\begin{align*}
g^+(z)+g^+(e^{\pi i}\,z)=\log(1-e^{-2\pi iz})\,;
\end{align*}
for any $z\in \cH=S(0,\pi)$,
\begin{align*}
g^+(z)+g^+(e^{-\pi i}\,z)=\log(1-e^{2\pi iz})\,.
\end{align*}
\end{lemma}

\begin{proof}
The result follows immediately from \eqref{equation:g+-} and \eqref{equation:g-+}.
\end{proof}

Lemma \ref{lemma:g} is essentially the Euler's reflection formula on $\Gamma$-function, as it is easy to see that, from \eqref{equation:IGamma}, $I_\Gamma(z)=-g^+(z)$. If we set $G(\tau,\nu)=g^+(\frac\nu\tau)$, that is to say:
\begin{equation}\label{equation:GGamma}
G(\tau,\nu)=-\log\Gamma(\frac\nu\tau+1)+\bigl(\frac\nu\tau+\frac12\bigr)\,\log\frac\nu\tau-\frac\nu\tau+\log\sqrt{2\pi}\,,
\end{equation}
then $G(\tau,\nu)$ is well defined and analytic over the domain $U^+$ given below:
\begin{equation}\label{equation:U}
U^+:=\{(\tau,\nu)\in\CC^*\times\CC^*: \nu/\tau\notin(-\infty,0]\}\,.
\end{equation}

\begin{proposition}\label{proposition:G} Let $G$ be as in \eqref{equation:GGamma}. Then, for any $(\tau,\nu)\in U^+$,
\begin{equation}\label{equation:G+}
G(\tau,\nu)+G(\tau,-\nu)=\log(1-e^{\mp 2\pi i \nu/\tau})
\end{equation}
according to $\displaystyle \frac\nu\tau\in S(-\pi,0)$ or $S(0,\pi)$, respectively.
\end{proposition}

\begin{proof}
The statement comes from Lemma \ref{lemma:g}.
\end{proof}

\subsection{Modular-type expansion of $(x;q)_\infty$}\label{subsec:complex}

We shall discuss how to understand Theorem \ref{theorem:mainter} in the complex plane, for both complex $q$ and complex $x$. As before, let $\tilde\CC^*$ be the Riemann surface of the logarithm function.
Let $\cM$ be the automorphism of the $2$-dimensional complex manifold $\tilde\CC^*\times\tilde\CC^*$ given as follows:
\begin{equation*}
\cM: (q,x)\mapsto\cM(q,x)=\bigl(\iota(q),\iota_q(x)\bigr),
\end{equation*}
where
\begin{equation}\label{equation:q*z*}
\iota(q)=q^*:=e^{4\pi^2/\log q},\qquad
\iota_q(x)=x^*:=e^{2\pi i\log x/\log q}\,.
\end{equation}
In the following, we will use the notations $q^*$ and $x^*$ instead of $\iota(q)$ and $\iota_q(x)$ each time when any confusion does not occur.

If we let $\tilde D^*=\exp\bigl(i\cH\bigr)\subset\tilde\CC^*$, then $\cM$ induces an automorphism over the sub-manifold $\tilde D^*\times \tilde\CC^*$. From then now, we always write $q=e^{2\pi i\tau}$, $x=e^{2\pi i\nu}$ and suppose $ \tau\in\cH$, so that $0<|q|<1$.
Sometimes we shall use the pairs of modular variables $(\tau^*,\nu^*)$ as follows:
\begin{equation}\label{equation:modularvariables}
i(\tau)=\tau^*:=-1/\tau,\quad
i_\tau(\nu)=\nu^*:=\nu/\tau\,,
\end{equation}
so that we can continue to write $q^*=e^{2\pi i\tau^*}$ and $x^*=e^{2\pi i\nu^*}$\,.

\begin{theorem}\label{theorem:xqmain}
Let $q=e^{2\pi i\tau}$, $x=e^{2\pi i\nu}$ and let $q^*$, $x^*$ as in \eqref{equation:q*z*}. The following relation holds for any $\tau\in\cH$ and $\nu\in\CC\setminus\bigl((-\infty,-1]\cup[1,\infty)\bigr)$ such that $\nu/\tau\notin(-\infty,0]$:
\begin{align}\label{equation:xqmain}
(x;q)_\infty=&q^{-1/24}\,\sqrt{1-x}\,\,(x^*q^*;q^*)_\infty\cr
&\times\,\exp\Bigl(\frac{\li(x)}{\log q}+G(\tau,\nu)+P(\tau,\nu)\Bigr)\,,\end{align}
where $\sqrt{1-x}$ stands for the principal branch of $e^{\frac12\log(1-x)}$, $\li$ denotes the dilogarithm recalled in \eqref{equation:dilog}, $G$ is given by \eqref{equation:GGamma} and where $P$ denotes the function $P_-$ defined in Lemma \ref{lemma:Pd}.
\end{theorem}

\begin{proof}
By Theorem \ref{theorem:MP} and relation \eqref{equation:PP-}, we arrive at the expression
$$
M(-\frac{\log q}{2\pi},\frac{\log x}{\log q})=\log(x^*\,q^*;q^*)_\infty+P(\tau,\nu)\,;
$$
making then suitable variable change in \eqref{equation:mainter} allows one to arrive at \eqref{equation:xqmain}, by taking into account the standard analytic continuation argument.
\end{proof}

If we denote by $G^*$ the anti-symetrization of $G$ given by
\begin{equation*}
G^*(\tau,\nu)=\frac12\,\bigl(G(\tau,\nu)-G(\tau,-\nu)\bigr)\,,
\end{equation*}
then, according to relation \eqref{equation:G+}, we may rewrite \eqref{equation:xqmain} as follows:
\begin{align}\label{equation:xqmainH}
(x;q)_\infty=&q^{-1/24}\,\sqrt{\frac{1-x}{1-x^*}}\,\,(x^*;q^*)_\infty\cr
&\times\,\exp\Bigl(\frac{\li(x)}{\log q}+G^*(\tau,\nu)+P(\tau,\nu)\Bigr)
\end{align}
if $\nu\in \tau\cH$, and
\begin{align}\label{equation:xqmainH-}
(x;q)_\infty=&q^{-1/24}\,\frac{\sqrt{(1-x)(1-1/x^*)}}{1-x^*}\,\,(x^*;q^*)_\infty\cr
&\times\,\exp\Bigl(\frac{\li(x)}{\log q}+G^*(\tau,\nu)+P(\tau,\nu)\Bigr)
\end{align}
if $\nu\in-\tau\cH$, that is, if $\displaystyle\frac\nu\tau\in S(-\pi,0)$.

In the above,  $G^*$ and $P$ are odd functions on the variable $\nu$:
\begin{equation}\label{equation:symetricG*P}
G^*(\tau,-\nu)=-G^*(\tau,\nu),\quad
P(\tau,-\nu)=-P(\tau,\nu)\,;
\end{equation}
$\li$ satisfies the so-called \emph{Landen's transformation} \cite[Theorem 2.6.1, p. 103]{AAR}:
\begin{equation}\label{equation:Landen}
\li(1-x)+\li(1-\frac1x)=-\frac12\,\bigl(\log x\bigr)^2\,.
\end{equation}

Finally, if we write $\vec\omega=(\omega_1,\omega_2)=(1,\tau)$ and denote by $\Gamma_2(z,\vec\omega)$ the Barnes' double Gamma function associated to the double period $\vec\omega$ (\cite{Ba}), then Thoerme \ref{theorem:xqmain} and Proposition 5 of \cite{Sh} imply that
\begin{align}\label{equation:doubleGamma}
\frac{\Gamma_2(1+\tau-\nu,\vec\omega)}{\Gamma_2(\nu,\vec\omega)}=&\sqrt{i}\,\sqrt{1-x}\,\exp\Bigl(\frac{\pi i}{12\tau}+\frac{\pi i}{2}\bigl(\frac{\nu^2}{\tau}-(1+\frac1\tau)\nu\bigr)\cr
&+\frac{\li(x)}{\log q}+G(\tau,\nu)+P(\tau,\nu)\Bigr)\cr
=&\,\sqrt{2\sin\pi\nu}\,\,\exp\Bigl(\frac{\pi i}{12\tau}+\frac{\nu(\nu-1)\pi i}{2\tau}\cr
&+\frac{\li(e^{2\pi i\nu})}{2\pi i\tau}+G(\tau,\nu)+P(\tau,\nu)\Bigr)\,. \end{align}

\subsection{Completed Ramanujan's formula of $(x;q)_\infty$}\label{subsection:Ramanujan} For any positive integer $n$, consider the following power series of $z$:
\begin{equation}\label{equation:An0}
A_n(z)=\sum_{k\ge 0}\frac{B_{2(n+k)}}{2(n+k)\,(2k+1)!}\,z^{2k+1}\,,
\end{equation}
where $B_{2m}$ denote the Bernoulli numbers. We recall the following identity \cite[p.12, (1.2.10)]{AAR}:
\begin{equation}\label{equation:Bernoullin}
B(\frac{t}{2\pi})=\sum_{m\ge 1}\frac{B_{2m}}{(2m)!}\,t^{2m}\,.
\end{equation}
where $B(t)$ is as given in \eqref{equation:B}. Moreover, from the relation \cite[p.29, (1.6.4)]{AAR}
$$
\int_0^\infty\frac{t^{2m}}{e^{2\pi t}-1}\,\frac{dt}{t}=(-1)^{m-1}\,\frac{B_{2m}}{4m}\,,
$$
one finds that
\begin{equation}\label{equation:An}
 A_n(z)=(-1)^{n-1}\,\int_0^\infty\frac{2t^{2n}}{e^{2\pi t}-1}\,\sin(tz)\,dt\,.
\end{equation}

\begin{proposition}\label{proposition:An}
For each positive integer $n$, the above-defined function $A_n$ can be continued into an analytic function for all $z\in\CC$ such that $z^2\notin(-\infty,-4\pi^2]\subset\RR$.
\end{proposition}

\begin{proof}
At the right hand side of \eqref{equation:An}, one can replace the integration loop $[0,\infty)$ with any half-line $[0,\infty e^{id})$ in the half-plane $\Re t>0$, with $d\in(-\frac\pi2,\frac\pi2)$. By noticing that this integral converges for all $z\in\CC$ verifying $\vert\Re (ize^{id})\vert<\Re(2\pi e^{id})$, one gets the expected analytic continuation domain.
\end{proof}

Let $\theta$ be the (formal) power series given as in \eqref{equation:Ramtheta}; in view of \eqref{equation:An0}, it follows that
\begin{equation}\label{equation:Ramtheta1}
\theta=\sum_{n\ge 1}\frac{B_{2n}\,A_n(\log x)}{(2n)!}\,(\log q)^{2n-1}\,,
\end{equation}
where $\log x\in\CC\setminus\Delta$, with $\Delta=[2\pi i,\infty i)\cup(-\infty i,-2\pi i]$. Remember the notation $S(a,b)$ is introduced in \eqref{equation:Sab}; by making use of the formula \eqref{equation:xqmain} of Theorem \ref{theorem:xqmain}, Ramanujan's formula \eqref{equation:Ram} may be {\it analytically} completed as follows.

\begin{theorem}\label{theorem:Ram}
Let $q$, $x$, $q^*$, $x^*$, $\tau$ and $\nu$ as in Theorem \ref{theorem:xqmain}.
If $P=P_-(\tau,\nu)$ be as given in Lemma \ref{lemma:Pd} and $s=\log x/\log q=\nu/\tau$, then it follows that
\begin{equation}\label{equation:Ram1}
 (x;q)_\infty=\frac{\sqrt{2\pi s(1-x)}}{\Gamma(s+1)}\,q^{-1/24}\,e^{s(\log s-1)}\,e^{\frac{\li(x)}{\log q}+P}\,(q^*x^*;q^*)_\infty\,.
\end{equation}
Moreover, the power series $\theta$ defined in \eqref{equation:Ramtheta} or \eqref{equation:Ramtheta1} is a uniform asymptotic expansion of $(-P)$ at $\log q=2\pi i\tau=0$ in the following sense: for any $\epsilon\in(0,\frac\pi2)$, there exists a positive constant $C_\epsilon>0$ such that, for any positive integer $N$, if
$$
K_N(\nu)=\int_0^{+\infty}\frac{\vert\hbox{\bf sh}(\nu t)\vert}{e^t-1}\,t^{2N+1}\,\frac{dt}{t}\,,
$$
then  the following estimates hold for all integer $N\ge 0$ and all $(\tau,\nu)\in S(\epsilon,\pi-\epsilon)\times\CC$ verifying $\vert\Re\nu\vert<1$:
\begin{equation}\label{equation:Ramasymp}
\bigl\vert P(\tau,\nu)+\sum_{n=1}^N \frac{B_{2n}\,A_n(2\pi i\nu)}{(2n)!}\,(2\pi i\tau)^{2n-1}\bigr\vert\le \frac{C_\epsilon\,K_N(\nu)}{(2\pi-\epsilon)^{2N}}\,\vert\tau\vert^{2N+1}.
\end{equation}
\end{theorem}

\begin{proof}
Relation \eqref{equation:Ram1} comes directly from \eqref{equation:xqmain} and \eqref{equation:GGamma}.

To prove \eqref{equation:Ramasymp}, we come back to the relation \eqref{equation:cotB} of the proof of Lemma \ref{lemma:Pd}, and let $f$ be  defined in $\CC\setminus (2\pi\ZZ)$ as follows:
$$
f(t)=\cot\frac t2-\frac2t=2i B\bigl(\frac{it}{2\pi}\bigr)\,.
$$
Obviously, $f$ may be analytically continued at $t=0$. By taking into account the relation \cite[p. 12, (1.2.10)]{AAR}, we find that, for $\vert t\vert<1$,
$$
B(t)=\sum_{n\ge 1}B_{2n}\,\frac{(2\pi t)^{2n-1}}{(2n)!}\,,
$$
which implies that, if $|t|<2\pi$, then
\begin{equation}\label{equation:expansionf}
f(t)=2\sum_{n\ge 1}(-1)^n\,\frac{B_{2n}}{(2n)!}\,t^{2n-1}\,.
\end{equation}
Moreover, the following limit holds for any given $\epsilon\in(0,\frac\pi2)$:
\begin{equation}\label{equation:limitf}
\lim_{S(-\pi+\epsilon,-\epsilon)\ni t\to\infty}f(t)=i\,. 
\end{equation}

Let $\epsilon\in(0,\frac\pi2)$. Let $f_0(t)=0$ and for any positive integer $N$, let
$$
f_N(t)=2\sum_{n=1}^N(-1)^n\,\frac{B_{2n}}{(2n)!}\,t^{2n-1}\,.
$$
By \eqref{equation:limitf} and the fact that $f$ is an odd analytic function for $\vert t\vert<2\pi$, one may find a positive constant $C_\epsilon>0$ such that, for all $t\in S(-\pi+\epsilon,-\epsilon)$ and all positive integer $N$:
\begin{equation}\label{equation:RamafN}
\bigl\vert f(t)-f_N(t)\bigr\vert\le C_\epsilon\,(2\pi-\epsilon)^{-2N}\,t^{2N+1}\,. 
\end{equation}
Choose any $d\in(-\pi+\epsilon.-\epsilon)$, consider the relation \eqref{equation:Pd} and let
$$P_N(\tau,\nu)=\int_0^{\infty e^{id}}\frac{\sin\frac\nu\tau t}{e^{it/\tau}-1}\,f_N(t)\,\frac{dt}t$$
and
$$
R_N(\tau,\nu)=P(\tau,\nu)-P_N(\tau,\nu)\,.
$$
By \eqref{equation:An}, one finds easily that
$$
P_N(\tau,\nu)=-\sum_{n=1}^N \frac{B_{2n}\,A_n(2\pi i\nu)}{(2n)!}\,(2\pi i\tau)^{2n-1}\,.
$$
To finish the proof, one needs only to give estimates for $R_N$, that may be easily done with the help of \eqref{equation:RamafN}. We omit the details.
\end{proof}

Remark finally the asymptotic expansion \eqref{equation:Ramasymp} is valid for all $(x;q)\in{\CC^*}^2$ such that $\vert q\vert<1$. In fact, if $x=e^{2\pi\nu}$, one may always suppose that $\Re\nu\in[0,1)$.

\subsection{Some remarks when $q$ tends toward one}\label{subsec:limit}

For the sake of simplicity, we will limit ourself to the real case and we suppose $q\to 1^-$ by real values in $(0,1)$, so that one can let $\tau=i\alpha$, $\alpha\to 0^+$. As $\tau^*=-1/\tau=i/\alpha$, one may observe that $\Im(\tau^*)\to+\infty$ and therefore $q^*\to0^+$ rapidly or, exactly saying, exponentially with respect to the variable $1/\alpha$. The relation
$$|x^*|=e^{2\pi\nu/\alpha}=e^{2\pi\,\Re(\nu)/\alpha}$$ shows that, as $\alpha\to0^+$, the modular variable $x^*$ belongs to the unit circle if and only if the initial variable $x$ takes a real value; otherwise,  $x^*$ goes rapidly to $\infty i$ or $0$ according to the sign of $\Re\nu$.
Let $\Re(\nu)\in[0,1)$; it follows that
$$
(x^*q^*;q^*)_\infty=1+O\bigl(e^{-(1-\Re(\nu))/\alpha}\,\bigr)\,\sim 1.
$$

Moreover, as $q\to 1-$, $q^*=e^{4\pi^2/\log q}$ becomes {\it exponentially small} and the expression $(x^*q^*;q^*)_\infty$ can not be represented by any semi-convergent power series of  $\log q$.  If we compare \cite[Entry 6, p. 265]{Be1} with \cite[Entry 6', p. 268]{Be1} in only which the equality symbol ``=" is used, we would  like to believe Ramanujan really wanted to give a convergent expression to $\log(x;q)_\infty$.

Finally, if $\tau=i\alpha$, $\alpha>0$, it is easy to verify that the following limits hold: for any fix $\nu\in(0,1)$,$$
\lim_{\alpha\to 0^+}P(\tau,\nu)=\lim_{\alpha\to 0^+}G(\tau,\nu)=0\,.
$$
Remark that the first one can be deduced immediately from \eqref{equation:Ramasymp} with $N=0$.
Therefore, by Theorem \ref{theorem:Ram}, we find:
$$
\log(x;q)_\infty\sim \frac{\log({1-x})}2\,-\frac{\li(x)}{2\pi\alpha}
$$
when $q=e^{-2\pi\alpha}\to 1^-$\,.

In a forthcoming paper, we shall give a compactly uniform Gevrey asymptotic expansion for $(x;q)_\infty$ when $q\to 1$ inside the unit disc, $x$ being a complex parameter; see \cite[\S 1.4.1, p. 84-86]{MR} for Gevrey asymptotic expansion with parameters.

\section{Dedekind $\eta$-function, Jacobi $\theta$-function and Lambert series}\label{section:modular}

In the following, we will first see in what manner Theorem \ref{theorem:xqmain} essentially contains the modular equations known for $\eta$ and $\theta$-functions; see Theorems \ref{theorem:eta} and \ref{theorem:theta}. In \S \ref{subsec:L12}, we will consider two families of series, called $L_1$ and $L_2$, that can be obtained as logarithmic derivatives of the infinite product $(x;q)_\infty$; some modular-type relations will be given in Theorem \ref{theorem:L12modular}. In \S \ref{subsec:Lambert}, classical Lambert series will be viewed as particular cases of the previous series $L_1$ and $L_2$.

\subsection{Dedekind $\eta$-function}\label{subsec:eta}

Let us mention a first application of Theorem \ref{theorem:xqmain} as follows.

\begin{theorem}\label{theorem:eta}
Let $q=e^{2\pi i\tau}$, $ \tau\in\cH$ et let $q^*=e^{-2\pi i/\tau}$. Then
\begin{equation}\label{equation:eta}
(q;q)_\infty=q^{-1/24}\,\sqrt{\frac i\tau}\,\,(q^*)^{1/24}\,(q^*;q^*)_\infty\,.
\end{equation}
\end{theorem}

\begin{proof}If we set
\begin{align*}
G_0(\tau,\nu)=&G(\tau,\nu)-\log\sqrt{\frac{2\pi\nu}\tau}\cr
=&\log\Gamma(\nu^*+1)+\nu^*\log\nu^*-\nu^*\,\,,
\end{align*}
we can write relation \eqref{equation:xqmain} of Theorem \ref{theorem:xqmain} as follows:
\begin{align*}
(xq;q)_\infty=&q^{-1/24}\,\sqrt{\frac{2\pi\nu}{(1-e^{2\pi i\nu})\tau}}\,\,(x^*q^*;q^*)_\infty\cr
&\times\,\exp\Bigl(\frac{\li(x)}{\log q}+G_0(\tau,\nu)+P(\tau,\nu)\Bigr)\,,
\end{align*}
where $x=e^{2\pi i\nu}$. Suppose $\nu\to 0$, so that $x\to 1$, $\nu^*\to 0$ and $x^*\to 1$; from \eqref{equation:GGamma}, it follows:
$$
\lim_{\nu\to 0}G_0(\tau,\nu)=0\,;
$$
therefore, one easily gets relation \eqref{equation:eta}, remembering that $\displaystyle \li(1)=\frac{\pi^2}6$ and that $P(\tau,0)=0$ as is said in \eqref{equation:symetricG*P}.
\end{proof}

The function $(q;q)_\infty$ plays a very important role in number theory and it is really linked with the well-known Dedekind $\eta$-function:
\begin{equation}\label{equation:etadefinition}
\eta(\tau)=e^{\pi\tau i/12}\prod_{n=1}^\infty(1-e^{2n\pi\tau i})=q^{1/24}\,(q;q)_\infty\,,
\end{equation}
where $\tau\in\cH$. For instance, see \cite[Lectures VI, VIII]{Ha} and \cite[Chapters 10, 11]{AAR}. The modular relation \eqref{equation:eta}, written as
$$
\eta(-\frac1\tau)=\sqrt{\frac{\tau}i}\,\,\eta(\tau)\,,
$$
 is traditionally obtained as consequence of Poisson's summation formula ({\it cf} \cite[p. 597-599]{Gu}) or that of Mellin transform of some Dirichlet series ({\it cf} \cite[p. 538-540]{AAR})\,; see also \cite{Si}, for a simple proof.

\subsection{Modular relation on Jacobi theta function}\label{subsec:theta}

In order to get the modular equation for Jacobi theta function, we first mention the following relation for any $x\in \CC\setminus(-\infty,0]$:
\begin{equation}\label{equation:dilog+}
\li(-x)+\li(-\frac1x)=-\frac{\pi^2}6-\frac12\,\bigl(\log x\bigr)^2\,,
\end{equation}
which can be deduced directly from the definition \eqref{equation:dilog} of $\li$, for
$$
\frac{d\ }{dx}\bigl(\li(-x)+\li(-\frac1x)\bigr)=\frac{\log(1+x)}{-x}-\frac{\log(1+\frac1x)}{-x}=\frac{\log x}{-x}\,.
$$

One can also check \eqref{equation:dilog+} by making use of a suitable variable change and considering both the Landen's transformation \eqref{equation:Landen} and formula \cite[(2.6.6), p. 104]{AAR}:
$$
\li(x)+\li(1-x)=\frac{\pi^2}6-\log x\,\log(1-x)\,;
$$
see \cite{Za} for more information.

For any fix $q=e^{2\pi i\tau}$, $\tau\in\cH$, the modular variable transformation $x\mapsto \iota_q(x)=x^*$ introduced in \eqref{equation:q*z*} defines an automorphism of the Riemann surface $\tilde\CC^*$ of the logarithm and  satisfies the following relations ($q^*=\iota(q)=e^{4\pi^2/\log q}$):
\begin{equation}\label{equation:modularvariables+}
\iota_q(xy)=\iota_q(x)\,\iota_q(y);\quad
\iota_q(e^{2k\pi i})=(q^*)^{-k}\,,\quad
\iota_q(q^k)=e^{2k\pi i}
\end{equation}
for all $k\in\RR$. In particular, one finds:
\begin{equation}\label{equation:modularvariables-1}
\iota_q(\sqrt q\,x)=e^{\pi i}\,\iota_q(x)\,,\quad
\iota_q(xe^{i\pi})=\iota_q(x)/\sqrt{q^*}\,.
\end{equation}

As usual, for any $m$ given complex numbers $a_1$, $\dots$, $a_m$, let
$$
(a_1,\ldots,a_m;q)_\infty=\prod_{k=1}^m(a_k;q)_\infty\,.
$$

\begin{theorem}\label{theorem:theta}
Let $q=e^{2\pi i\tau}$ and $x=e^{2\pi i\nu}$ and let
\begin{equation}\label{equation:thetadefinition}
\theta(q,x)=(q,-\sqrt q\, x,-\frac{\sqrt q}x;q)_\infty\,.
\end{equation}
Then, the following relation holds for any $\tau\in\cH$ and any $\nu$ of the  Riemann surface of the logarithm:
\begin{equation}\label{equation:thetamodular}
\theta(q,x)=q^{1/8}\,\sqrt{\frac i{\tau x}}\,\,\exp\Bigl(-\frac{(\log \frac x{\sqrt q}\,)^2}{2\log q}\,\Bigr)\,\theta(q^*,x^*)\,.
\end{equation}

\end{theorem}

\begin{proof} First, suppose $\nu\in \tau\cH$ and write $(x;q)_\infty$ and $(1/x;q)_\infty$ by means of \eqref{equation:xqmainH} and \eqref{equation:xqmainH-}, respectively. By taking into account relation \eqref{equation:symetricG*P} about the parity of $G^*$ and $P$,  we  find:
\begin{align*}
(x,\frac1x;q)_\infty=&q^{-1/12}\,\frac{1-x}{1-1/x^*}\,\,\sqrt{-\frac1x}\,\,(x^*,\frac1{x^*};q^*)_\infty\,\cr
&\,\exp\Bigl(\frac1{\log q}\,\bigl(\li(x)+\li(\frac1x)\bigr)\Bigr)\,,
\end{align*}
where we used the relation $1/x^*=(1/x)^*$, deduced from \eqref{equation:modularvariables+}. Thus, it follows that
$$
(xq,\frac1x;q)_\infty=q^{-1/12}\,\sqrt{-\frac1x}\,\,(x^*,\frac{q^*}{x^*};q^*)_\infty\,\exp\Bigl(\frac1{\log q}\,\bigl(\bigl(\li(x)+\li(\frac1x)\bigr)\Bigr)\,,
$$
Change $x$ by $e^{-\pi i}\,x/\sqrt q$ and make use of \eqref{equation:modularvariables-1} and \eqref{equation:dilog+}; we get:
\begin{align*}
(-\sqrt q\, x,-\frac{\sqrt q}x;q)_\infty=&q^{-1/12}\,\sqrt{\frac{\sqrt q}x\,}\,\,(-\sqrt{q^*}\,x^*,-\frac{\sqrt{q^*}}{x^*};q^*)_\infty\cr
&\times\,\exp\Bigl(\frac1{\log q}\,\bigl(-\frac{\pi^2}6-\frac12\,(\log\frac x{\sqrt q} \,)^2\bigr)\Bigr)
\end{align*}
for $x\notin(-\infty,0]$.
By the modular equation \eqref{equation:eta} for $\eta$-function, we arrive at the expected formula \eqref{equation:thetamodular}.

Finally we end the proof of the Theorem by the standard analytic continuation argument.
\end{proof}

Known as the modular formula on Jacobi's theta function, relation \eqref{equation:thetamodular} can be written as follows:
$$
\theta(q,x)=\sqrt{\frac i{\tau }}\,\,\exp\Bigl(-\frac{(\log x\,)^2}{2\log q}\,\Bigr)\,\theta(q^*,x^*)\,,
$$
which has a very long history, and is attached to Gauss, Jacobi, Dedekind, Hermite, {\it etc} .... It is generally obtained by applying Poisson's summation formula to the Laurent series expansion:
$$
(q,-\sqrt q\, x,-\frac{\sqrt q}x;q)_\infty=\sum_{n\in\ZZ}q^{\frac{n^2}2}\,x^n\,,
$$
which is the so-called Jacobi's triple product formula; for instance, see \cite[\S 10.4, p. 496-501]{AAR}.

\subsection{Another proof of Theta modular equation}\label{subsec:anotherproof}

As what is pointed out in \cite[p. 214-215]{zh2}, formula \eqref{equation:thetamodular} can be interpreted in term of $q$-difference equations. We shall elaborate on this idea and give a simple proof for \eqref{equation:thetamodular}.

For any fix $q$ such that $0<|q|<1$, let
$$f_1(x)=\theta(q,x), \quad
f_2(x)=\displaystyle\sqrt\frac1{x}\,\exp\Bigl(-\frac{(\log \frac x{\sqrt q}\,)^2}{2\log q}\,\Bigr)
$$
and
$$
f(x)=\displaystyle \frac{f_1(x)}{f_2(x)}=g(x^*)\,,
$$
where $x^*$ is given by \eqref{equation:q*z*}.
As $f_1$ and $f_2$ are solution of the same first order linear equation
$$y(qx)=\displaystyle\frac{1}{\sqrt q\,x}\,y(x),
$$
$f$ is a $q$-constant, that means, $f(qx)=f(x)$\,; equivalently, $g$ is uniform on the variable $x^*$, for $qx$ is translated into $x^* e^{2\pi i}$ by \eqref{equation:modularvariables+}. On the other hand, it is easy to check the following relation:
$$f(xe^{-2\pi i})=e^{-2\pi i\frac{\log x}{\log q}\,-\frac{2\pi^2}{\log q}}\,f(x)=\frac{1}{\sqrt {q^*}\,x^*}\,f(x),
$$
so that, using \eqref{equation:modularvariables+}, we find:
$$
g(q^*\,x^*)=\frac{1}{\sqrt {q^*}\,x^*}\,g(x^*)\,.
$$

Summarizing, $g$ is a uniform solution of $y(q^*x^*)=y(x^*)/({\sqrt {q^*}\,x^*})$ and vanishes over the $q^*$-spiral $-\sqrt{q^*}\,{q^*}^\ZZ$ of the $x^*$-Riemann surface of the logarithm; it follows that there exists a constant $C$ such that $g(x^*)=C\,\theta(q^*,x^*)$ for all $x^*\in\tilde\CC^*$. Write
$$
C=\frac{\theta(q,x)}{\theta(q^*,x^*)}\,\sqrt{x}\,\exp\Bigl(\frac{(\log\frac{x}{\sqrt{q}})^2}{2\log q}\Bigr)
$$
and let $x\to e^{\pi i}/\sqrt q$, so that $x^*\to e^{\pi i}/\sqrt{q^*}$. Since
$$
\frac{\theta(q,x)}{1+\sqrt q\, x}\to (q;q)^3_\infty\,,\quad
\frac{\theta(q^*,x^*)}{1+\sqrt{q^*}\,x^*}\to (q^*;q^*)_\infty^3
$$
and
$$
\lim_{x\to e^{\pi i}/\sqrt{q}}\frac{1+\sqrt q\, x}{1+\sqrt{q^*}\, x^*}=\frac{\sqrt{q}}{\sqrt{q^*}\,}\,\bigl(\frac{dx^*}{dx}\bigl\vert_{x=e^{\pi i}/\sqrt q}\bigr)^{-1}=\frac{\tau}{i}\,,
$$
where $q=e^{2\pi i\tau}$, we get the following expression, deduced from $\eta$-modular equation \eqref{equation:eta}:
\begin{align*}
C=q^{1/4}\,e^{-\frac{\pi i}{2}-\frac{\pi^2}{2\log q}}\,\frac{(q;q)_\infty^3}{(q^*;q^*)_\infty^3}\,\frac{\tau}{i}=q^{1/8}\,\sqrt\frac{i}{\tau}\,\,,
\end{align*}
and we end the proof of \eqref{equation:thetamodular}.

One key point of the previous proof is to use the dual variables $q^*$ and $x^*$; the underlying idea is really linked with the concept of local monodromy group of linear $q$-difference equations \cite[\S 2.2.3, Th\'eor\`eme 2.2.3.5] {Sa}. In fact, as there exists two generators  for the fundamental group of the elliptic curve $\CC^*/q^\ZZ$, one needs to consider the ``monodromy operators'' in two directions or ``two periods'', $x\mapsto x e^{2\pi i}$ and $x\mapsto xq$, which exactly correspond to $x^*\mapsto x^*q^*$ and $x^*\mapsto x^* e^{-2\pi i}$, in view of \eqref{equation:modularvariables+}.

\subsection{Generalized Lambert series $L_1$ and $L_2$}\label{subsec:L12}As before, let $q=e^{2\pi i\tau}$, $x=e^{2\pi i\nu}$ and suppose $\tau\in\cH$. Consider the following series, which can be considered as generalized Lambert series:
\begin{equation}\label{equation:L12}
L_1(\tau,\nu)=\sum_{n=0}^\infty \frac{x\,q^n}{1-x\,q^n}\,,\quad
L_2(\tau,\nu)=\sum_{n=0}^\infty \frac{(n+1)x\,q^n}{1-x\,q^n}\,,
\end{equation}
that are both absolutely convergent for any $x\in\CC\setminus q^{-\NN}$, due to the fact $|q|<1$.
By expanding each term $(1-xq^n)^{-1}$ into geometric series, one easily finds:
\begin{equation}\label{equation:transformL12}
L_1(\tau,\nu)=\sum_{n=0}^\infty\frac{x^{n+1}}{1-q^{n+1}}\,,\quad
L_2(\tau,\nu)=\sum_{n=0}^\infty\frac{x^{n+1}}{(1-q^{n+1})^2}\,,
\end{equation}
where convergence requires $x$ to be inside the unit circle $\vert x\vert<1$ of $x$-plane.

Observe that
\begin{equation}\label{equation:L12difference1}
 L_1(\tau,\nu+\tau)-L_1(\tau,\nu)=-\frac{x}{1-x}
\end{equation}
and
\begin{equation}\label{equation:L12difference2}
L_2(\tau,\nu+\tau)-L_2(\tau,\nu)=-L_1(\tau,\nu)\,.
\end{equation}
In this way, one may guess how to define more series such as $L_3$, $L_4$, {\it etc} ...

A direct computation yields the following formulas:
\begin{equation}\label{equation:L12'1}
L_1(\tau,\nu)=-x\,\frac{\partial\ }{\partial x}\log(x;q)_\infty\,,
\end{equation}
\begin{equation}\label{equation:L12'2}
L_2(\tau,\nu)=-q\,\frac{\partial\ }{\partial q}\log(x;q)_\infty+L_1(\tau,\nu)\,;
\end{equation}
\begin{align}\label{equation:xderivatives}
x\,\frac{\partial x^*}{\partial x}=\frac{x^*}\tau,\quad
x\,\frac{\partial \nu}{\partial x}=\frac1{2\pi i}
\end{align}
and
\begin{align}\label{equation:qderivatives}
q\,\frac{\partial x^*}{\partial q}=-\frac{\nu}{\tau^2}\,x^*,\quad
q\,\frac{\partial q^*}{\partial q}=\frac{1}{\tau^2}\,q^*,\quad
q\,\frac{\partial \tau}{\partial q}=\frac{1}{2\pi i}\,.
\end{align}
Here and in the following, $q$ and $x$ are considered as independent variables as well as the pair $(\tau,\nu)$ or their modular versions $(q^*,x^*)$ and $(\tau^*,\nu^*)$.

\begin{theorem}\label{theorem:L12modular} Let $q=e^{2\pi i\tau}$, $x=e^{2\pi i\nu}$ and let $q^*$, $x^*$, $\tau^*$ and $\nu^*$ be as in \eqref{equation:q*z*} and \eqref{equation:modularvariables}. If $\tau\in\cH$ and $\nu/\tau\not\in(-\infty,0]$, then the following relations hold:
\begin{align}\label{equation:L1modular}
L_1(\tau,\nu+\tau)=&\frac{\log(1-e^{2\pi i\nu})}{2\pi\tau i}+\frac {e^{2\pi i\nu}}{2(1-e^{2\pi i\nu})}+L_1(\tau^*,\nu^*+\tau^*)\,\frac1\tau\cr
&+\frac{1}{2\pi\tau i}\,\Bigl(\frac{\Gamma'(\frac\nu\tau+1)}{\Gamma(\frac\nu\tau+1)}-\log\frac\nu\tau-\frac\tau{2\nu}-\tau\,\frac{\partial\ }{\partial\nu}P(\tau,\nu)\Bigr)
\end{align}
and\begin{align}\label{equation:L2modular}L_2(\tau,\nu+\tau)=&\frac1{24}-\frac{\li(e^{2\pi i\nu})}{4\pi^2}\,\frac1{\tau^2}-L_1(\tau^*,\nu^*+\tau^*)\,\frac{\nu}{\tau^2}\cr
&+L_2(\tau^*,\nu^*+\tau^*)\,\frac1{\tau^2}-\frac{\nu}{2\pi  i\tau^2}\,\Bigl(\frac{\Gamma'(\frac\nu\tau+1)}{\Gamma(\frac\nu\tau+1)}\cr
&-\log\frac\nu\tau-\frac\tau{2\nu}+\frac{\tau^2}\nu\,\frac{\partial\ }{\partial\tau}P(\tau,\nu)\Bigr)\,.
\end{align}
\end{theorem}

\begin{proof} By taking the logarithmic derivative with respect to the variable $x$ in \eqref{equation:xqmain} and in view of \eqref{equation:xderivatives}, we find:
\begin{align*}
x\frac{\partial\ }{\partial x}\log(x;q)_\infty=&-\frac x{2(1-x)}+\frac{x^*}{1-x^*}\,\frac{1}{\tau}+\frac1\tau\,x^*\frac{\partial\ }{\partial x^*}\log(x^*;q^*)_\infty\cr
&-\log_q(1-x)+\frac{1}{2\pi i}\,\bigl(\frac{\partial\ }{\partial\nu}G(\tau,\nu)+\frac{\partial\ }{\partial\nu}P(\tau,\nu)\bigr)\,,
\end{align*}
so that, by  \eqref{equation:L12'1}, we arrive at the following expression:
\begin{align*}
L_1(\tau,\nu)=&\log_q(1-x)+\frac x{2(1-x)}+\bigl(L_1(\tau^*,\nu^*)-\frac{x^*}{1-x^*}\bigr)\,\frac1\tau\cr
&-\frac{1}{2\pi i}\,\bigl(\frac{\partial\ }{\partial\nu}G(\tau,\nu)+\frac{\partial\ }{\partial\nu}P(\tau,\nu)\bigr)\,.
\end{align*}
From \eqref{equation:GGamma}, it follows that
\begin{equation}\label{equation:G'nu}
\tau\,\frac{\partial\ }{\partial\nu}G(\tau,\nu)=-\frac{\Gamma'(\frac\nu\tau+1)}{\Gamma(\frac\nu\tau+1)}+\log\frac\nu\tau+\frac\tau{2\nu},
\end{equation}
that leads to the wanted relation \eqref{equation:L1modular}, using \eqref{equation:L12difference1}.

On the other hand, using \eqref{equation:qderivatives}, a direct computation shows that \eqref{equation:xqmain} implies the following expression:
\begin{align*}
q\frac{\partial\ }{\partial q}\log(x;q)_\infty=&-\frac1{24}+\frac\nu{\tau^2}\,x^*\frac{\partial\ }{\partial x^*}\,\log(1-x^*)-\frac\nu{\tau^2}\,x^*\frac{\partial\ }{\partial x^*}\,\log(x^*;q^*)_\infty\cr
&+\frac1{\tau^2}\,q^*\frac{\partial\ }{\partial q^*}\,\log(x^*;q^*)_\infty-\frac{\li(x)}{(\log q)^2}\cr
&+\frac1{2\pi i}\,\Bigl(\frac{\partial\ }{\partial\tau}G(\tau,\nu)+\frac{\partial\ }{\partial\tau}P(\tau,\nu)\Bigr),
\end{align*}
or, by taking into account \eqref{equation:L12'2}:
\begin{align}\label{equation:qqxq}
q\frac{\partial\ }{\partial q}\log(x;q)_\infty=&-\frac1{24}+\frac{\li(x)}{4\pi^2}\,\frac1{\tau^2}-\frac{x^*}{1-x^*}\,\frac{\nu}{\tau^2}+L_1(\tau^*,\nu^*)\,\frac{\nu+1}{\tau^2}\cr
&-L_2(\tau^*,\nu^*)\,\frac1{\tau^2}+\frac1{2\pi i}\,\bigl(\frac{\partial\ }{\partial\tau}G(\tau,\nu)+\frac{\partial\ }{\partial\tau}P(\tau,\nu)\bigr)\,.
\end{align}
Since
$$
\tau \frac{\partial\ }{\partial\tau}G(\tau,\nu)=-\nu \frac{\partial\ }{\partial\nu}G(\tau,\nu),
$$
by putting together \eqref{equation:L12'2}, \eqref{equation:L1modular}, \eqref{equation:G'nu} with the last relation \eqref{equation:qqxq}, we get the following relation:
\begin{align*}
L_2(\tau,\nu)=&\frac1{24}-\frac{\li(x)}{4\pi^2}\,\frac1{\tau^2}+\frac{x^*}{1-x^*}\,\frac{\nu-\tau}{\tau^2}-L_1(\tau^*,\nu^*)\,\frac{\nu+1}{\tau^2}\cr
&+L_1(\tau,\nu)+L_2(\tau^*,\nu^*)\,\frac1{\tau^2}-\frac{\nu}{2\pi  i\tau^2}\,\Bigl(\frac{\Gamma'(\frac\nu\tau+1)}{\Gamma(\frac\nu\tau+1)}\cr
&-\log\frac\nu\tau-\frac\tau{2\nu}+\frac{\tau^2}\nu\,\frac{\partial\ }{\partial\tau}P(\tau,\nu)\Bigr) \,,
\end{align*}
so that we arrive at the wanted formula \eqref{equation:L2modular} by applying \eqref{equation:L12difference2}.
\end{proof}

\subsection{Classical Lambert series}\label{subsec:Lambert}

If we let $\nu=\tau$,  we reduce series $L_1$ and $L_2$ to the following  classical Lambert series:
$$
L_1(\tau,\tau)=\sum_{n=0}^\infty\frac{q^{n+1}}{1-q^{n+1}}
$$
and
$$
L_2(\tau,\tau)=\sum_{n=0}^\infty\frac{(n+1)q^{n+1}}{1-q^{n+1}}=\sum_{n=0}^\infty\frac{q^{n+1}}{(1-q^{n+1})^2}\,.
$$
By considering the limit case $\nu=0$ in \eqref{equation:L1modular} and \eqref{equation:L2modular} of Theorem \ref{theorem:L12modular}, we arrive at the following statement.

\begin{theorem}\label{theorem:L120}
For all $\tau\in\cH$, the following relations hold:
\begin{align}\label{equation:L10}
L_1(\tau,\tau)=&\frac{\log(-2\pi i\tau)}{2\pi i\tau}+\frac{1}{4}\cr
&-\frac{\gamma}{2\pi\,i\tau}-\frac1{2\pi i}\,\frac{\partial\ }{\partial\nu}P(\tau,0)+L_1(\tau^*,\tau^*)\,\frac1\tau\,,
\end{align}
where $\gamma$ denotes Euler's constant, and
\begin{align}\label{equation:L20}
L_2(\tau,\tau)=&\frac{1}{24}+\frac1{4\pi i\tau}-\frac1{24\,\tau^2}+L_2(\tau^*,\tau^*)\,\frac1{\tau^2}\,.
\end{align}
\end{theorem}

\begin{proof}
If we set
$$
A(\tau,\nu)=\frac{1}{2\pi i\tau}\,\log\frac{1-e^{2\pi i\nu}}{\nu/\tau}-\frac12\,\bigl(\frac{e^{2\pi i\nu}}{1-e^{2\pi i\nu}}+\frac{1}{2\pi i\nu}\bigr)\,,
$$
we can write \eqref{equation:L1modular} as follows:
\begin{align*}
L_1(\tau,\nu+\tau)=&A(\tau,\nu)+L_1(\tau^*,\nu^*+\tau^*)\,\frac1\tau\cr&+\frac{1}{2\pi\tau i}\,\Bigl(\frac{\Gamma'(\frac\nu\tau+1)}{\Gamma(\frac\nu\tau+1)}-\tau\,\frac{\partial\ }{\partial\nu}P(\tau,\nu)\Bigr)\,,
\end{align*}
so that, remembering $\gamma=-\displaystyle\Gamma'(1)$, we get \eqref{equation:L10}, as it is easy to see that
$$
\lim_{\nu\to 0}A(\tau,\nu)=\frac{\log(-2\pi i\tau)}{2\pi i\tau}+\frac{1}{4}\,.
$$

In the same time, putting $\nu=0$ in \eqref{equation:L2modular} allows one to obtain relation \eqref{equation:L20}, for $P(\tau,0)=0$ for all $\tau\in\cH$ implies $\displaystyle\frac{\partial\ }{\partial\tau}P(\tau,0)=0$ identically.
\end{proof}

Formula \eqref{equation:L10} has been known since Schl\"omilch; see Stieltjes \cite[(84), p. 54]{St}. Relation \eqref{equation:L20} is really a modular relation and is traditionally obtained by taking derivative with respect to the variable $\tau$ in modular formula \eqref{equation:eta}; see \cite[Exercises 6 and 7, p. 71]{Ap}.

\section{Proof of Theorem \ref{theorem:main}}\label{section:proof}

 In all this section, we let
$$
q=e^{-a}=e^{-2\pi\alpha},\quad
x=e^{-(1+\xi)a}=q^{1+\xi}
$$
and suppose
$$
a=2\pi\alpha>0,\quad
0<q<1,\quad
\xi> -1,\quad
0<x<1\,,
$$
For any positive integer $N$, define
\begin{equation}\label{equation:VN}
V_N(a,\xi):=\sum_{n=1}^N\log(1-e^{-(n+\xi)a})\,.
\end{equation}
It is esay to see that
$$
\log\,(x;q)_\infty=\lim_{N\to\infty}V_N(a,\xi)\,.
$$

We shall prove Theorem \ref{theorem:main} in several steps, and our approach is well inspired by Stieltjes' work \emph{\'Etude de la fonction $P(a)=\displaystyle\sum_{1}^\infty\frac1{e^{\frac na}-1}$} that one can find in his Thesis \cite[p. 57-62]{St}. The starting point is to use the fact that $\displaystyle\frac1{e^{\sqrt{2\pi }\,u}-1}-\frac1{\sqrt{2\pi}\,u}$ is a self-reciprocal function with respect to Fourier sine transform  \cite[(7.2.2), p.~179]{Ti}, so that one may write each finite sum $V_N$ by four or five appropriate sine or cosine integrals depending of $N$ and make then estimation over these integrals.

\subsection{Some preparatory formulas}\label{subsec:preparatory}

We are going to use the following formulas:
\begin{equation}\label{equation:sinus}
\int_0^\infty\frac{\sin\lambda u}{e^{2\pi u}-1}du=\frac14+\frac{1}2\bigl(\frac1{e^\lambda-1}-\frac1\lambda\bigr)\,,
\end{equation}
and
\begin{equation}\label{equation:cosinus}
\int_0^\infty\frac{1-\cos\lambda u}{e^{2\pi u}-1}\frac{du}u=\frac\lambda4+\frac{1}{2}\log\frac{1-e^{-\lambda}}\lambda\,,
\end{equation}where $\lambda$ is assumed to be a real or complex number such that $\vert\Im \lambda\vert<2\pi$; notice that the second formula can be deduced from the first one by integrating on $\lambda$. For instance, see  \cite[(82) \& (83), p. 57]{St}, \cite[(7.2.2), p.~179]{Ti} or \cite[Example 2, p. 122]{WW}.

From \eqref{equation:cosinus}, it follows that
\begin{align}\label{equation:RN}
V_N(a,\xi)=&\sum_{n=1}^N\log\,(n+\xi)-\frac N2\, \xi a\cr
&-\frac{N(N+1)}4\, a+N\log\, a+R_N(a,\xi)\,,
\end{align}
where
$$
R_N(a,\xi)=\int_0^\infty\frac{h_N(au,\xi)}{e^{2\pi u}-1}\,\frac{du}u
$$
and
$$
h_N(t,\xi)=2N-2\sum_{n=1}^N\cos(n+\xi)t\,.
$$
By using the elementary relations
$$
2\sum_{n=1}^N\cos nt=\cos Nt+\sin Nt\cot\frac t2-1
$$
and
$$
2\sum_{n=1}^N\sin nt=\sin Nt +(1-\cos Nt)\cot\frac t2\,,
$$
we obtain:
$$
h_N(t,\xi)=2N+\cos\xi t-\cos( N+\xi)t+\bigl(\sin \xi t-\sin(N+\xi)t\bigr)\cot\frac t2\,.
$$
Let us define the following integrals:
\begin{equation}\label{equation:RN1}
 R_{N}^{(1)}(a,\xi)=\int_0^\infty\frac{\cos\xi au-\cos( N+\xi)au}{e^{2\pi u}-1}\,\frac{du}u
\end{equation}
and
\begin{equation}\label{equation:RN2}
 R_{N}^{(2)}(a,\xi)=\int_0^\infty\frac{2N+\bigl(\sin \xi au-\sin(N+\xi)au\bigr)\cot\frac {au}2}{e^{2\pi u}-1}\,\frac{du}u\,,
\end{equation}
so that
$$
R_N(a,\xi)=R_{N}^{(1)}(a,\xi)+R_{N}^{(2)}(a,\xi).
$$

We will look for the limits of $R_N^{(1)}$ and $R_N^{(2)}$ while $N$ becomes indefinitely large. To simplify, we will write
$
a_N\sim_N b_N
$ if the quantity $(a_N-b_N)$ tends to zero as $N\to\infty$.
From \eqref{equation:cosinus}, we first observe the following result.

\begin{lemma}\label{lemma:RN1}
The following relation holds:
\begin{equation}\label{equation:limitRN1}
R_N^{(1)}(a,\xi)\sim_N\frac N4\,a-\frac12\,\log N-\frac12\,\log\frac{1-e^{-\xi a}}\xi\,.
\end{equation}
\end{lemma}

\begin{proof}
Applying \eqref{equation:cosinus} to the following integrals
$$
\int_0^\infty\frac{1-\cos\xi au}{e^{2\pi u}-1}\,\frac{du}u\,,\quad
\int_0^\infty\frac{1-\cos (N+\xi) au}{e^{2\pi u}-1}\,\frac{du}u
$$
implies directly Lemma \ref{lemma:RN1}.
\end{proof}

The following well-known result, due to Riemann, will be often taken into account in the course of the proof.

\begin{lemma}\label{lemma:Riemann}
Let $f$ be a continuous and integrable function on a finite or infinite closed interval $[\alpha,\beta]\subset\overline{\RR}$. Then the following relations hold:
$$
\int_\alpha^\beta f(t)\sin Nt\,dt\sim_N0,\quad
\int_\alpha^\beta f(t)\cos Nt\,dt\sim_N0.
$$
\end{lemma}

\subsection{First part of $R_N^{(2)}$} \label{subsec:RN21}

The integral \eqref{equation:RN2} of $R_N^{(2)}$ seems more complicated than $R_N^{(1)}$, because of the simple poles at $u=\displaystyle\frac2a\,\pi$, $\displaystyle\frac4a\,\pi$, $\displaystyle\frac6a\,\pi$, {\it etc}, $\cdots$ that the function $\displaystyle\cot\frac{au}2$ admits on $(0,+\infty)$. In such a situation, one very classical technique may consist in replacing the function by its decomposition in simple parts as given in \eqref{equation:cot}.
By considering $\displaystyle\frac2{au}$ instead of $\displaystyle\cot\frac{au}2$ in \eqref{equation:RN2}, we are led to the following integral:
\begin{equation}\label{equation:RN21}
R_N^{(21)}(a,\xi)=\frac2a\,\int_0^\infty\frac{Nau+\sin \xi au-\sin(N+\xi)au}{e^{2\pi u}-1}\,\frac{du}{u^2}\,;
\end{equation}
if we set
\begin{equation}\label{equation:RN22}
R_N^{(22)}(a,\xi)=4a\,\int_0^\infty\sum_{n=1}^\infty\frac{\sin(N+\xi)au-\sin \xi au}{4\pi^2n^2-a^2u^2}\,\frac{du}{e^{2\pi u}-1}\,,
\end{equation}
then, in view of \eqref{equation:cot}  we obtain the following equality:
\begin{equation*}
R_N^{(2)}(a,\xi)=R_N^{(21)}(a,\xi)+R_N^{(22)}(a,\xi)\,.
\end{equation*}

\begin{lemma}\label{lemma:RN21}
The following relation holds for $a>0$ and $\xi>0$:
\begin{align}
\label{equation:limitRN21}
R_N^{(21)}(a,\xi)\sim_N&\frac{N(N+2\xi)}4\,a-(N+\xi)\log(N+\xi)+N(1-\log a)\cr
&-\frac{\pi^2}{6a}+\xi\log\xi-\frac1a\int_0^{\xi a}\log(1-e^{-t})\,dt\,.
\end{align}
\end{lemma}

\begin{proof}
For any pair $(N,\xi)\in\NN\times\RR$, let $f_N(a,\xi)=\displaystyle\frac a2R_N^{(21)}(a,\xi)$; it is easy to see that $a\mapsto f_N(a,\xi)$ represents an odd analytic function at the origin $0$ of the real axis, for merely
$$
f_N(a,\xi)=\int_0^\infty\frac{Nau+\sin \xi au-\sin(N+\xi)au}{e^{2\pi u}-1}\,\frac{du}{u^2}\,.
$$
Let $f_N'(a,x)$ denote the derivative of $f_N(a,\xi)$ with respect to the variable $a$. It follows that
$$
f_N'(a,\xi)=\int_0^\infty\frac{(N+\xi)\bigl(1-\cos(N+\xi)au\bigr)-\xi(1-\cos \xi au)}{e^{2\pi u}-1}\,\frac{du}{u}\,,
$$
so that applying \eqref{equation:cosinus} gives rise to the following relation:
\begin{align*}
f_N'(a,\xi)=&\frac{N(N+2\xi)}4\,a+\frac{N+\xi}2\,\log\frac{1-e^{-(N+\xi)a}}{N+\xi}\cr
&-\frac N2\,\log a-\frac\xi2\,\log\frac{1-e^{-\xi a}}{\xi}\,.
\end{align*}

To come back to $f_N(a,x)$, we integrate $f_N'(t,\xi)$ over the interval $(0,a)$ and remark that $f_N(0,\xi)=0$; it follows that
\begin{align}\label{equation:fa}
f_N(a,\xi)=&\frac{N(N+2\xi)}8\,a^2-\frac{N+\xi}2\,a\,\log(N+\xi)-\frac N2\,(\log a-1)a \cr
&+\frac\xi2\,(\log\xi)\,a+\frac12\,I(a,N+\xi)-\frac12\,I(a,\xi)\,,
\end{align}
where
$$
I(a,\delta)=\int_0^{\delta a}\log(1-e^{-t})\,dt\,.
$$

Now we suppose  $a>0$ and let $N\to+\infty$. Noticing that
\begin{equation}\label{equation:limitI}
I(a,N+\xi)\sim_N\int_0^{\infty}\log(1-e^{-t})\,dt=-\li(1)=-\frac{\pi^2}6,
\end{equation}
we get immediately \eqref{equation:limitRN21} from \eqref{equation:fa}.

\end{proof}

The term $\displaystyle-\frac{\pi^2}{6a}$ included in expression \eqref{equation:RN21} plays a most important role for understanding the asymptotic behavior of $\log (x;q)_\infty$ as $q\to 1^-$, that is, $a\to 0^+$. The crucial point is formula \eqref{equation:limitI}, that remains valid for all complex numbers $a$ such that $\Re a>0$.

\subsection{Intermediate part in $R_N^{(2)}$} \label{subsec:RN22}

Now consider $R_N^{(22)}(a,\xi)$ of \eqref{equation:RN22}; then
\begin{equation}\label{equation:RN22IJ}
R_N^{(22)}(a,\xi)=\frac2\pi\, \bigl(I_N(a,\xi)-J_N(a,\xi)\bigr)\,,
\end{equation}
if we set
\begin{equation}\label{equation:IN}
I_N(a,\xi)=\int_0^\infty\sum_{n=1}^\infty\frac{\sin 2nN\pi t\cos2n\xi\pi t}{n(e^{4n\pi^2t/a}-1)}\,\frac{dt}{1-t^2}
\end{equation}
and
\begin{equation}\label{equation:JN}
J_N(a,\xi)=\int_0^\infty\sum_{n=1}^\infty\frac{\sin2n\xi\pi t(1-\cos2nN\pi t)}{n(e^{4n\pi^2t/a}-1)}\,\frac{dt}{1-t^2}\,,
\end{equation}
Here, each series under $\displaystyle \int$ converges absolutely to an integrable function over $(0,\infty)$ excepted maybe  near zero and $t=1$. Lemma \ref{lemma:ANM} given below will tell us how to regularize the situation at origin; notice also that these integrals behave \emph{more convergent} at $t=1$ than $0$, due to big factors $(e^{4n\pi^2t/a}-1)$.

\begin{lemma}\label{lemma:ANM}
Let $\delta\in(0,1)$, $\lambda>0$ and let $\{h_{n,N}\}_{n,N\in\NN}$ be a uniformly bounded family of continuous functions on $[0,\delta]$. For any positive integer $M$, let $A_M({N})$ denote the integral given by
$$
A_{M}(N)=\int_0^\delta\sum_{n=1}^M\frac{h_{n,N}(t)}{n}\,\bigl(\frac1{e^{n\lambda t}-1}-\frac 1{n\lambda t}\bigr)\,\frac{dt}{1-t^2}\,.
$$
Then, as $M\to \infty$, the sequence $\{A_M(N)\}$  converges uniformly for $N\in\NN$.
\end{lemma}

\begin{proof}
We suppose $\lambda=1$, the general case being analogous; thus, one can write $A_{M}(N)$ as follows:
$$
A_{M}(N)=\sum_{n=1}^M\int_0^{n\delta}{h_{n,N}(t/n)}\,\bigl(\frac1{e^{t}-1}-\frac 1{t}\bigr)\,\frac{dt}{n^2-t^2}\,
$$
Observe that the function $(e^t-1)^{-1}-t^{-1}$ increases rapidly from $-1/2$ toward zero when $t$ tends to infinity by positive values; indeed, $(e^t-1)^{-1}-t^{-1}=O(t^{-1})$ for $t\to+\infty$. Therefore, if we make use of the relation $\displaystyle \int_0^{n\delta}=\int_0^{\sqrt n\,\delta}+\int_{\sqrt n\,\delta}^{n\delta}$ and let $n\to+\infty$, we find:
$$
\Bigl\vert\int_0^{n\delta}{h_{n,N}(t/n)}\,\bigl(\frac1{e^{t}-1}-\frac 1{t}\bigr)\,\frac{dt}{n^2-t^2}\Bigr\vert\le C\,n^{-3/2}\,,
$$
where $C$ denotes a suitable positive constant independent of $N$ and $n$; this ends the proof of Lemma \ref{lemma:ANM}.
\end{proof}

We come back to the integral $I_N$ given in \eqref{equation:IN}.
\begin{lemma}\label{lemma:IN}
The following relation holds:
\begin{equation}\label{equation:limitIN}
I_N(a,\xi)\sim_N\frac \pi{48}\,a-\frac\pi2\sum_{n=1}^\infty\frac{\cos2n\pi\xi}{n(e^{4n\pi^2/a}-1)}\,.
\end{equation}
\end{lemma}

\begin{proof}
We fix a small $\delta>0$, cut off the interval $(0,\infty)$ into four parts $(0,\delta)$, $(\delta,1-\delta)$, $(1-\delta,1+\delta)$ and $(1+\delta,\infty)$, and the corresponding integrals will be denoted by $I_N^{0\delta}$, $I_N^{\delta+}$, $I_N^{1\mp\delta}$ and $I_N^{\delta\infty}$, respectively. According to Lemma \ref{lemma:Riemann}, we find:
\begin{equation}\label{equation:limitIN24}
I_N^{\delta+}(a,\xi)\sim_N0,\qquad
	I_N^{\delta\infty}(a,\xi)\sim_N0,
\end{equation}
	for
	$$
	\sum_{n=M}^\infty\bigl(\int_\delta^{1-\delta}+\int_{1+\delta}^\infty\bigr)\frac{|\sin 2nN\pi t\cos2n\xi\pi t|}{n(e^{4n\pi^2t/a}-1)}\,\frac{dt}{|1-t^2|}\to 0
	$$
when $M\to\infty$.	In the same way, we may observe that
\begin{align}\label{equation:IN1}
I_N^{0\delta}(a,\xi)&\sim_N \frac a{4\pi^2}\,\int_0^\delta\sum_{n=1}^\infty\frac{\sin 2nN\pi t\cos 2n\xi\pi t}{n^2t({1-t^2})}\,{dt}\cr
	&\sim_N \frac a{4\pi^2}\,\int_0^\delta\sum_{n=1}^\infty\frac{\sin 2nN\pi t}{n^2}\,\frac{dt}t\cr
&=\frac a{4\pi^2}\,\sum_{n=1}^\infty\,\int_0^\delta\frac{\sin 2nN\pi t}{n^2}\,\frac{dt}t\,,
	\end{align}
	where the first approximation relation is essentially obtained  from Lemma \ref{lemma:ANM}, combining together with Lemma \ref{lemma:Riemann}.
	Since
\begin{equation}\label{equation:sintt}
	\int_0^\infty\frac{\sin t}t\,dt=\frac\pi2,
\end{equation}
from \eqref{equation:IN1} we deduce the following relation:
\begin{align}\label{equation:limitIN1}
I_N^{0\delta}(a,\xi)\sim_N \frac a{8\pi}\sum_{n=1}^\infty\frac1{n^2}=\frac{\pi}{48}\,a\,.
	\end{align}

A similar analysis allows one to write the following relations for the remaining integral $I_N^{1\mp\delta}$:
	\begin{align}\label{equation:limitIN3}
	I_N^{1\mp\delta}(a,\xi)&\sim_N\frac12\, \int_{1-\delta}^{1+\delta}\sum_{n=1}^\infty\frac{\sin 2nN\pi t\cos2n\xi\pi t}{n(e^{4n\pi^2t/a}-1)}\,\frac{dt}{1-t}\cr
&	\sim_N\frac12\,\int_{1-\delta}^{1+\delta}\sum_{n=1}^\infty\frac{\cos2n\xi\pi}{n(e^{4n\pi^2/a}-1)}\,\frac{\sin 2nN\pi t}{1-t}\,dt\cr
&	\sim_N\frac12\,\sum_{n=1}^\infty\frac{\cos2n\xi\pi}{n(e^{4n\pi^2/a}-1)}\,\int_{-\infty}^{+\infty}-\frac{\sin t}{t}\,dt\cr
&=-\frac\pi2\,\sum_{n=1}^\infty\frac{\cos2n\xi\pi}{n(e^{4n\pi^2/a}-1)}\,,
\end{align}
where the last equality comes from \eqref{equation:sintt}.

Accordingly, we obtain the wanted expression \eqref{equation:limitIN} by putting together the estimates \eqref{equation:limitIN24}, \eqref{equation:limitIN1} and \eqref{equation:limitIN3} and thus the proof is complete.\end{proof}

\subsection{Singular integral as limit part of $R_N^{(2)}$} \label{subsec:RN23}

In order to give estimates for $J_N(a,\xi)$ of \eqref{equation:JN}, we shall make use of the Cauchy principal value of a singular integral. The situation we have to consider is the following \cite[\S 6.23, p. 117]{WW}:  $f$ be a continuous function over $(0,1)\cup(1,+\infty)$ such that, for any $\epsilon>0$, $f$ is integrable over both intervals $(0,1-\epsilon)$ and $(1+\epsilon,+\infty)$; one defines
$$
{\cal PV}\int_0^\infty f(t)dt=\lim_{\epsilon\to 0^+}\bigl(\int_0^{1-\epsilon}f(t)\,dt+\int_{1+\epsilon}^\infty f(t)\,dt\bigr)
$$
whenever the last limit exists.

\begin{lemma}\label{lemma:JN}
The following relation holds:
\begin{equation}\label{equation:limitJN}
J_N(a,\xi)\sim_N{\cal PV}\int_0^\infty\sum_{n=1}^\infty\frac{\sin2n\xi \pi t}{n(e^{4n\pi^2t/a}-1)}\,\frac{dt}{1-t^2}\,.
\end{equation}
\end{lemma}

\begin{proof}
For any given number $\epsilon\in(0,1)$, let
$$
J_N^{(1\mp\epsilon)}(a,\xi)=\int_{1-\epsilon}^{1+\epsilon}\sum_{n=1}^\infty\frac{\sin2n\xi \pi t(1-\cos2nN\pi t)}{n(e^{4n\pi^2t/a}-1)}\,\frac{dt}{1-t^2}\,;
$$
Thanks to suitable variable change, we can get the following expression:
$$
J_N^{(1\mp\epsilon)}(a,\xi)=\int_0^\epsilon\sum_{n=1}^\infty\bigl(\frac{h_n(a,\xi,t)}{2-t}-\frac{h_n(a,\xi,-t)}{2+t}\bigr)\,\frac{1-\cos2nN\pi t}{nt}\,dt\,,
$$
where
$$
h_n(a,\xi,t)=\frac{\sin2n\xi \pi(1 -t)}{e^{4n\pi^2(1-t)/a}-1}\,.
$$
From Lemma \ref{lemma:Riemann}, one deduces:
\begin{equation}\label{equation:JN1}
J_N^{(1\mp\epsilon)}(a,\xi)\sim_N
\int_0^\epsilon\sum_{n=1}^\infty\bigl(\frac{h_n(a,\xi,t)}{2-t}-\frac{h_n(a,\xi,-t)}{2+t}\bigr)\,\frac{dt}{nt}\,.
\end{equation}

Again applying Lemma \ref{lemma:Riemann} implies that
\begin{align*}
J_N(a,\xi)-J_N^{(1\mp\epsilon)}(a,\xi)\sim_N\bigl(\int_0^{1-\epsilon}+\int_{1+\epsilon}^\infty\bigr)\sum_{n=1}^\infty\frac{\sin2n\xi \pi t}{n(e^{4n\pi^2t/a}-1)}\,\frac{dt}{1-t^2}\,,
\end{align*}
which, using \eqref{equation:JN1}, permits to conclude, as it is clear that
$$
\lim_{\epsilon\to0^+}\int_0^\epsilon\sum_{n=1}^\infty\bigl(\frac{h_n(a,\xi,t)}{2-t}-\frac{h_n(a,\xi,-t)}{2+t}\bigr)\,\frac{dt}{nt}=0\,.
$$
\end{proof}

\subsection{End of the proof of Theorem \ref{theorem:main}}

\begin{proof}
Consider the functions $V_N(a,\xi)$ given in \eqref{equation:VN} and recall that $\log(x;q)_\infty$ is the limit of $V_N(a,\xi)$ when $N$ goes to infinity; so we need to know the limit behavior of the right hand side of \eqref{equation:RN}  for infinitely large $N$.

Letting
$$
G_N(a,\xi)=\sum_{n=1}^N\log\,(n+\xi)-\frac N2\, \xi a-\frac{N(N+1)}4\, a+N\log\, a\,,
$$
it follows that
$$
V_N(a,\xi)=G_N(a,\xi)+R_N^{(1)}(a,\xi)+R_N^{(21)}(a,\xi)+\frac2\pi\,(I_N(a,\xi)-J_N(a,\xi))\,,
$$
where $R_N^{(1)}$, $R_N^{(21)}$, $I_N$ and $J_N$ are considered in Lemmas \ref{lemma:RN1}, \ref{lemma:RN21}, \ref{lemma:IN} and \ref{lemma:JN}, respectively. From Stirling's asymptotic formula \cite[Theorem 1.4.1, page 18]{AAR}, one easily gets:
\begin{align*}
\sum_{n=1}^N\log\,(n+\xi)=&\log\Gamma(N+\xi+1)-\log\Gamma(\xi+1)\cr
\sim_N&\log\sqrt{2\pi}+(N+\xi+\frac12)\log(N+\xi+1)\cr
&-(N+\xi+1)-\log\Gamma(\xi+1)\,.
\end{align*}
Thanks to \eqref{equation:limitRN21} of Lemma \ref{lemma:RN21}, one finds:
\begin{align*}
G_N(a,\xi)+R_N^{(21)}(a,\xi)\sim_N&-\frac N4\,a+\frac12\log N-\log\Gamma(\xi+1)
-\frac{\pi^2}{6a}\cr
&+\log\sqrt{2\pi}-\xi+\xi\log \xi-\frac1{a}\int_0^{\xi a}\log(1-e^{-t})\,dt\,.
\end{align*}
Thus, using \eqref{equation:RN1} of Lemma \ref{lemma:RN1}, one can deduce the following expression:
\begin{align*}
&G_N(a,\xi)+R_N^{(1)}(a,\xi)+R_N^{(21)}(a,\xi)\cr\sim_N&
-\frac{\pi^2}{6a}+\log\sqrt{2\pi}-\xi-\log\Gamma(\xi+1)\cr
&-\bigl(\xi+\frac12\bigr)\log\frac{1-e^{-\xi a}}\xi+\frac1{a}\int_0^{\xi a}\frac t{e^t-1}\,dt\,,
\end{align*}
which implies the starting formula \eqref{equation:main} of our paper with the help of Lemmas \ref{lemma:IN} and \ref{lemma:JN}, replacing all  $a$ by $2\pi\alpha$.
\end{proof}

{\bf Acknowledgements}

{The Author would like to express thanks to his friends and colleagues Anne Duval and Jacques Sauloy for their numerous valuable suggestions and remarks.
}

\end{document}